\documentclass[12pt, a4paper]{article}

\usepackage{amsmath}
\usepackage{amssymb}
\usepackage{mathrsfs}
\usepackage{epsfig}
\usepackage{xcolor}
\usepackage{graphicx}
\usepackage{subcaption}
\usepackage{enumerate,mathtools}
\usepackage{mathscinet}
\usepackage{cite}

\newcommand{\samethought}{\smallskip\noindent}
\newcommand{\newthought}{\medskip\noindent}
\newcommand{\qed}{\nolinebreak\hfill {$\Box$} \par\medbreak}
\newcommand{\finis}{\nolinebreak\hfill {$\bigtriangleup$}
  \par\medbreak}

\newcommand{\inprod}[2]{\langle #1 \mid #2 \rangle}

\newcommand{\eps}{\varepsilon}

\newcommand{\halfje}{{\textstyle \frac{1}{2}}}

\newcommand{\re}{{\mathop{\rm Re}}\,}
\newcommand{\im}{{\mathop{\rm Im}}\,}

\newcommand{\diag}{\mathop{\mathrm{diag}}}
\newcommand{\dd}[2]
   {\frac{\partial #1}{\partial #2}}

\newcommand{\E}{{\rm e}}
\newcommand{\I}{{\rm i}}

\newcommand{\R}{\mathbb{R}}
\newcommand{\N}{\mathbb{N}}
\newcommand{\Z}{\mathbb{Z}}

\newcommand{\C}{\mathbb{C}}
\newcommand{\T}{\mathbb{T}}

\newcommand{\KAM}{\textsc{kam}}

\newcommand{\Het}{\mathit{Het}}

\newcommand{\cO}{\mathcal O}

\newcommand{\sB}{\mathscr B}
\newcommand{\sP}{\mathscr P}

\newtheorem{theorem}{Theorem}[section]

\newtheorem{rmrk}[theorem]{Remark}
\newtheorem{xmpl}[theorem]{Example}

\newenvironment{remarks}{%
  \removelastskip \vspace*{9pt} \noindent
  \textbf{Remarks.}\begin{enumerate}}
  {\end{enumerate}}

\begin{document}

\title{\protect\Large Normal-normal resonances in a double Hopf
  bifurcation}

\author{ {\protect\normalsize Henk Broer}
  \protect\\[-1mm]
  {\protect\footnotesize\protect\it Bernoulli Institute for
    Mathematics, Computer Science and Artificial Intelligence}
  \protect\\[-2mm]
  {\protect\footnotesize\protect\it Rijksuniversiteit Groningen,
    9747~AG~Groningen, The Netherlands}
  \protect\\
  {\protect\normalsize
    Heinz Han{\ss}mann} \protect\\[-1mm]
  {\protect\footnotesize\protect\it Mathematisch Instituut,
    Universiteit Utrecht} \protect\\[-2mm]
  {\protect\footnotesize\protect\it Postbus 80010,
    3508~TA Utrecht, The Netherlands} \protect\\
  {\protect\normalsize Florian Wagener} \protect\\[-1mm]
  {\protect\footnotesize\protect\it Center for Nonlinear
    Dynamics in Economics and Finance (CeNDEF)} \protect\\[-2mm]
  {\protect\footnotesize\protect\it Amsterdam School of
    Economics, Universiteit van Amsterdam} \protect\\[-2mm]
  {\protect\footnotesize\protect\it Postbus 15867, 1001~NJ
    Amsterdam, The Netherlands}}

\date{\protect\normalsize 21 January 2020}

\maketitle

\begin{abstract}
  \noindent
  We investigate the stability loss of invariant
  $n$--dimensional quasi-periodic tori during a double Hopf
  bifurcation, where at bifurcation the two normal frequencies
  are in normal-normal resonance.
  Invariants are used to analyse the normal form approximations
  in a unified manner.
  The corresponding dynamics form a skeleton for the dynamics
  of the original system.
  Here both normal hyperbolicity and \KAM~theory are being used.
\end{abstract}

\section{Introduction}
\label{introduction}

Double Hopf bifurcations of $n$--dimensional tori in their simplest
form occur in $n+4$~dimensions.
Therefore we consider a family of vector fields on an $(n+4)$--dimensional
manifold that leave invariant a family of conditionally periodic
$n$--tori. 
Locally around the tori the dynamics are then of the form
\begin{subequations}
\label{notyetfloquet}
\begin{align}
   \dot{x} & \;\; = \;\; \omega(\mu) \; + \; \cO(z),
\label{notyetfloquetx}\\
   \dot{z} & \;\; = \;\; \Omega(\mu) \, z \; + \; \cO(z^2),
\label{notyetfloquetz}
\end{align}
\end{subequations}
where $\mu \in \R^s$ denotes the parameter.
We assume the tori to be reducible to Floquet form
and that the double Hopf bifurcation occurs at $\mu = 0$. 
The eigenvalues of $\Omega = \Omega(\mu)$ are the Floquet
exponents $\beta_j(\mu) \pm \I \alpha_j(\mu)$,
$\alpha_j, \beta_j \in \R$ and $j = 1, 2$, of the invariant
torus $\T^n \times \{ 0 \}$. 

\samethought
For $\mu = 0$ the tori are at a double Hopf singularity if
$\beta_j(0) = 0$ and $\alpha_j(0) \neq 0$ for $j = 1, 2$.
Within this singularity, the normal frequencies $\alpha_j(0)$
are at normal-normal resonance if moreover
\begin{equation}
  \label{eq:normal:resonance}
  \ell_2 \alpha_1(0) \;\; = \;\; \ell_1 \alpha_2(0)
\end{equation}
for some $0 \neq \ell = (\ell_1, \ell_2) \in \Z^2$.
The aim of this article is to study the dynamics of the
family~\eqref{notyetfloquet} near such resonances.

\newthought
The normal linear part of~\eqref{notyetfloquet} is
obtained by truncating the $\cO(z)$ and $\cO(z^2)$ terms in
\eqref{notyetfloquetx} and~\eqref{notyetfloquetz}, respectively.
The flow of this normal linear part is equivariant with respect to a
$\T^n \times \T^2$ action and in the resonant case 
\eqref{eq:normal:resonance} the symmetry group is even
larger (the symmetry group is maximal in the case
$\alpha_1 (0) = \alpha_2 (0) = 0$, which we do not consider further).
The aim of normal form theory is to push the symmetry to higher
order terms by successive transformations, compare
with~\cite{GH83}. 
In the non-resonant case the entire $\T^n \times \T^2$ symmetry
is inherited by the normalised dynamics. 
However, in the present normal-normal resonant case the normalised
terms only inherit a $\T^n \times \T^1$ symmetry.

\samethought
We apply finitely many normalising transformations so to obtain
a $\T^n \times \T^1$ symmetric truncation.
In particular, the coefficients in this truncation are
$x$--independent.
This decouples the normal dynamics from the internal dynamics. 
The normal dynamics retains a $\T$~symmetry.
After reducing this out, a three-dimensional system remains.
Equilibria of this reduced system correspond to either
$n$--dimensional or $(n+1)$--dimensional invariant tori in the
normal form truncation.
The full system is a small perturbation of this
truncation; to assess the impact of the higher order remainder terms,
we invoke both normal hyperbolicity and \KAM~theory.

\newthought
In~\eqref{eq:normal:resonance} we may assume that 
$0 < \ell_1 \le \ell_2$ and $\alpha_j(0) > 0$, $j = 1, 2$,
if necessary by relabeling the
frequencies and scaling the variables. 
We restrict to weak resonances: these are resonances
where in~\eqref{eq:normal:resonance} we have $0 < \ell_1 < \ell_2$.
This rules out the strongly resonant case $\ell = (1, 1)$.

\newthought
\begin{remarks}
\item[-] Instead of starting in $n+4$~dimensions we could also take as
starting point a family of vector fields on an $(n+m)$--dimensional
manifold.
However, our first step would then be to split the normal spectrum
of the bifurcating $n$--torus into the $4$~purely imaginary
eigenvalues $\pm \I \alpha_1(0), \pm \I \alpha_2(0)$ and the
remaining $m-4$~hyperbolic eigenvalues and to then perform a
restriction to the resulting $(n+4)$--dimensional normally
hyperbolic invariant manifold on which we have the situation
sketched above.
We note that this procedure generically leads to a finitely
differentiable system. 

\item[-] We generally assume our systems to be real analytic,
  which facilitates the application of \KAM~theory.
  However, for applications on centre manifolds finitely
  differentiable versions of \KAM~theory are
  available~\cite{BHW20,BHT90}.

\item[-] One of the improvements of normalisation theory in the analytic
  setting is that the remainder term can be made exponentially
  small~\cite{BHW18,BHJVW03} in the perturbation parameter. 
  In the finitely differentiable case this would be only polynomially small.
\end{remarks}
  
\subsection{Previous work }

The resonant double Hopf bifurcation has already been
studied by many authors.
The strong $1{:}1$~resonance, to which we do not contribute,
has been studied in~\cite{SD82, GKL90, NDLE94}.
Of the weak resonances, the 
$1{:}2$~resonance has gotten the most attention~\cite{KP88,
  LL96, KRS03, LPE03, V05, RAM13}; next to it, the 
$1{:}3$~resonance~\cite{LPE03} and the
$2{:}3$~resonance~\cite{RAM10, RAM15} have been
investigated as well.

\newthought
Most previous work has been on the case $n = 0$ where it is an
equilibrium that has two (normal) frequencies and we therefore
also speak of $1$--dimensional invariant tori instead of periodic
orbits.
The non-resonant case where only $\beta_1(0) = \beta_2(0) = 0$,
but $\alpha_1(0)$ and~$\alpha_2(0)$ do not
satisfy~\eqref{eq:normal:resonance} for any $0 \neq \ell \in
\Z^2$, has a $\T^2$~symmetry instead of a $\T^1$~symmetry in the normal
forms.
Reduction leads to a $2$--dimensional basis system and the
dynamics of this $2$--dimensional normal form was already laid out
in~\cite{GH83}.
However, it took until~\cite{li16} for a formal proof of persistence
of the resulting invariant tori under perturbation from the normal
form back to the original system.

\newthought
To summarise the literature on the weakly resonant case (again mostly
$n = 0$), it is convenient to call `resonance droplet' the part of
parameter space for which $1$--dimensional invariant tori with
non-zero amplitudes exist.
This corresponds for instance with the `mixed modes' situation
of Knobloch and Proctor~\cite{KP88}.
They made an extensive local study
of that part of the boundary of the resonance droplet of the
$1{:}2$ resonant double Hopf bifurcation where a $1$--dimensional
torus bifurcates to a $2$--dimensional torus, checking all
configurations and determining global bifurcations.
LeBlanc and Langford~\cite{LL96} used for the same Hopf
bifurcation Lyapunov-Schmidt theory and singularity theory to study
only these $1$--dimensional tori.
Luongo et al.~\cite{LPE03} applied a multiple timescale analysis
to the resonant $1{:}2$ and $1{:}3$ double Hopf bifurcation in a
straightforward perturbative approach.
They perform an extensive numerical analysis, finding several
instances of resonance droplets.   
Volkov~\cite{V05} takes $n \geq 2$ and studies the persistence of
$n$--dimensional quasi-periodic tori at a $1{:}2$ double Hopf normal
resonance.
Finally, Revel et al.\ \cite{RAM10, RAM13, RAM15} study $1{:}2$ and
$2{:}3$ resonant double Hopf bifurcations numerically. 

\newthought
\begin{remarks}
  
\item[-] Already in the work of Guckenheimer and Holmes~\cite{GH83} on the
  non-resonant case the treatment was devided into $12$~cases.
  In this paper we focus on the case labeled VIa in~\cite{GH83}, 
where a subordinate Hopf bifurcation takes place after reduction
of the $\T^n \times \T^2$~symmetry.
Next to this examplary case we note that our approach
also applies to all other cases as treated
by  Knobloch and Proctor~\cite{KP88}.

\item[-] Next to the subordinate bifurcations of co-dimension~$1$ we
  also meet bifurcations of co-dimension~$2$ that act as organising
  centers of the dynamics.
  These are subordinate resonant Hopf bifurcations, fold-Hopf
  bifurcations, and blue sky bifurcations.
  Also compare with~\cite{BDV07,BNRSW07}.
  The fold-Hopf bifurcation does not seem to have been noticed earlier
  in the literature and a full treatment is outside the present scope,
  so we refer to~\cite{BHW20}.
 
\item[-] In the sequel we introduce the invariants $\tau_j$, $j =
  1,2,3,4$.
  Here $\tau_1$ and $\tau_2$ generate the $\T^2$--symmetry, while
  $\tau_3$ and $\tau_4$ are a global version of the resonant angle.
  This systematic approach is standard for Hamiltonian systems,
  and has already been used in~\cite{KP88} for the $1{:}2$~resonance.
  We shall encounter $3$--dimensional versions of the Arnol'd resonance
  tongues as met in families of circle diffeomorphisms.
  The intersection with a transverse plane has the form of a droplet,
  see figures \ref{fig:SN1} and~\ref{fig:SN2} below.

\end{remarks}  

\subsection{Outline}

Our contribution is to establish the geometry of the resonance
droplet in a generic three-parameter unfolding of a general
$\ell_1{:}\ell_2$ resonant double Hopf bifurcation, excluding only
the strongly resonant $1{:}1$ situation.
We also put all occurring subordinate Hopf bifurcations on an equal
footing, but leave the occurring quasi-periodic fold-Hopf and heteroclinic
bifurcations to future research.

\samethought
To sketch our results we already refer to the bifurcation diagram
in figure~\ref{fig:hopf}, which describes the truncated
normal form after reduction of a $\T^n \times \T^2$~symmetry.
From this the symmetric skeleton is reconstructed.
The final step is to use perturbation theory back to the original
system.
This involves both normal hyperbolicity and \KAM~theory.

\section{Linear dynamics}
\label{ssec:linear:weak}

Consider on the phase space~$\R^4$ the system
\begin{align}\label{eq:linear}
  \dot{z} & \;\; = \;\; \Omega(\mu) z.
\end{align}
This is the normal part of the normal linear dynamics at the
invariant $n$--torus $\T^n \times \{ 0 \}$.

\samethought
We assume that the $4$ eigenvalues of $\Omega(\mu)$ are of the form
$\beta_j(\mu) \pm \I \alpha_j(\mu)$ and that they 
satisfy $\beta_j(0) = 0$ and  $\alpha_j(0) > 0$, for $j = 1, 2$.
Analogously to the internal frequency vector~$\omega(\mu)$, the
normal frequency vector of the invariant torus $\T^n \times \{ 0 \}$
is $\alpha(\mu) = (\alpha_1(\mu), \alpha_2(\mu))$.
Moreover, we assume that the normal frequencies are in
normal-normal resonance~\eqref{eq:normal:resonance}.
If necessary after rescaling time, it can be achieved that
$\alpha_1(0)$ takes some given value, for instance
$\alpha_1(0) = \ell_1$ (whence $\alpha_2(0) = \ell_2$).
Finally, we can assume that $\gcd(\ell_1, \ell_2) = 1$.
Normal-normal resonances at Hopf bifurcation generically
may occur in families depending on three or more parameters.

\newthought
Introducing complex variables $Z_j = z_{2j-1} + \I z_{2j}$,
$j=1,2$, at $\mu=0$ the system~\eqref{eq:linear} takes the form
\begin{equation*}
  \dot{Z}_j \;\; = \;\; \I \alpha_j(0) Z_j.
\end{equation*}
Under the linear flow there are four basic invariants
\begin{subequations}
\label{invariants}
\begin{alignat}{3}
  \label{invariants1}
  \tau_1 & \; = \; \halfje Z_1 \bar{Z}_1, \quad &
  \tau_3 & \; = \; \frac{\re Z_1^{\ell_2} \bar{Z}_2^{\ell_1}}
  {\ell_1! \ell_2!}, 
  \\
  \label{invariants2}
  \tau_2 & \; = \; \halfje Z_2 \bar{Z}_2, & 
  \tau_4 & \; = \; \frac{\im Z_1^{\ell_2} \bar{Z}_2^{\ell_1}}
  {\ell_1! \ell_2!},
\end{alignat}
\end{subequations}
i.e.\ all other functions that are invariant under the flow
of~\eqref{eq:linear} can be expressed as functions of the~$\tau_k$. 
These invariants are related by the syzygy
\begin{equation*}
  \tau_1^{\ell_2} \tau_2^{\ell_1}
  \, - \,
  G_\ell(\tau_3^2 + \tau_4^2) \; = \; 0,
  \quad
  G_\ell
  \; = \;
  \frac {(\ell_1!)^2 (\ell_2!)^2}{2^{\ell_1 + \ell_2}}.
\end{equation*}
The first two invariants are related to polar co-ordinates 
by the formula $Z_j = \sqrt{2 \tau_j} \E^{\I \phi_j}$.
At $\mu = 0$ the angles~$\phi_j$ have the equations of motion
\begin{displaymath}
   \dot{\phi}_j \;\; = \;\; \alpha_j (0), \quad j = 1, 2.
\end{displaymath}
As the frequencies are in resonance, the orbits
foliate each~$\T^2$ into closed $1$--dimensional orbits.
This is brought out explicitly by introducing resonance-adapted
angles $(\theta, \vartheta)$ such that 
\begin{subequations}
\label{eq:unimodular}
  \begin{align}
    \theta & \;\; = \;\; \ell_2 \phi_1 \; - \; \ell_1 \phi_2,
\label{eq:unimodular:theta} \\
    \vartheta & \;\; = \;\; m_1 \phi_1 \; + \; m_2 \phi_2.
\label{eq:unimodular:vartheta}
  \end{align}
\end{subequations}
While the right hand side of~\eqref{eq:unimodular:vartheta} is
clearly the inner product $\langle m \mid \phi \rangle$, we
define $\ell^{\perp} := (\ell_2, -\ell_1)$ to write the right
hand side of~\eqref{eq:unimodular:theta} as
$\langle \ell^{\perp} \mid \phi \rangle$

\samethought
In order that the equations \eqref{eq:unimodular} provide a
torus diffeomorphism, the $m_j$ have to be chosen such that the
matrix $\left(^{\ell_2}_{m_1} {}^{-\ell_1}_{m_2} \right)$ is
unimodular.
This is possible as $\gcd(\ell_1, \ell_2) = 1$.
Note that 
\begin{equation}
   \label{invariantphases}
   \tau_3 \; = \;
   \sqrt{G_\ell \tau_1^{\ell_2} \tau_2^{\ell_1}} \cos \theta 
   \quad \mbox{and} \quad
   \tau_4 \; = \;
   \sqrt{G_\ell \tau_1^{\ell_2} \tau_2^{\ell_1}} \sin \theta\, .
\end{equation}
In these variables, the linear dynamics are given by
\begin{equation}
\label{resonantinvariant}
\begin{array}{rcl}
  \dot{\tau}_j & = & 0, \quad j = 1, \ldots, 4 \\
  \dot{\theta} & = & 0 \quad \mbox{(redundant)} \\
  \dot{\vartheta} & = & m_1 \alpha_1(0) \; + \; m_2 \alpha_2(0)
               \;\; = \;\; 1,
\end{array}
\end{equation}
making $\vartheta$ a `fast' angle.
In appropriate co-ordinates the linear family
$\dot{z} = \Omega(\mu)z$ is given by
\begin{equation*}
  \Omega (\mu)\;\; = \;\;
  \begin{pmatrix}
    \beta_1(\mu) & -\alpha_1 (\mu) & 0 & 0 \\
    \alpha_1(\mu) & \beta_1 (\mu) & 0 & 0 \\
    0 & 0 & \beta_2 (\mu) & -\alpha_2 (\mu) \\
    0 & 0 &  \alpha_2 (\mu) & \beta_2 (\mu) \\
  \end{pmatrix}
  \, .
\end{equation*}
In terms of the complex variables this simplifies, first to
\begin{equation*}
  \dot{Z}_j \;\; = \;\;
  (\beta_j  (\mu) \, + \, \I \alpha_j (\mu)) Z_j,
  \quad j = 1, 2
\end{equation*}
and then to 
\begin{alignat*}{1}
  \dot{Z}_j & \;\; = \;\;
  (\I \alpha_j(0) \, + \, \beta_j + \I\delta_j) Z_j\, ,
\end{alignat*}
where the $\delta_j = \alpha_j (\mu) - \alpha_j (0)$ detune the
normal frequencies and form, together with
$\beta_j = \beta_j (\mu)$, $j = 1, 2$, the new independent
parameters.  
In fact, we have passed to
$\mu = (\delta_1, \delta_2, \beta_1, \beta_2, \nu)$,
where $\nu \in \R^{s-4}$ is a mute parameter. 
According to~\cite{arn71} this is a versal unfolding.

\section{Normal form}
\label{ssec:weak:nf}

We now add the higher order terms to the linear system.
On the phase space $\T^n \times \R^4$ the system becomes
\begin{alignat*}{3}
  \dot{x} & \;\; = \;\; f(x, z, \mu) & & \;\; = \;\;
  \omega(\mu) \; + \; \tilde f(x, z, \mu), \\
  \dot{z} & \;\; = \;\; h(x, z, \mu) & & \;\; = \;\;
  \Omega(\mu) z \; + \; \tilde h(x, z, \mu),
\end{alignat*}
where $\tilde f = \cO(|z|)$ and $\tilde h = \cO(|z|^2)$.
To this system are associated the vector field
\begin{equation}
  \label{eq:X}
  X \;\; = \;\; f \partial_x \; + \; h \partial_z
\end{equation}
and its normal linear part
$L = \omega \partial_z + \Omega z \partial_z$.
From the previous section we retain the assumptions on the
eigenvalues of $\Omega(\mu)$ and the complex variables~$Z_j$.
In these terms, the system takes the form
\begin{alignat*}{3}
  \dot{x} & \;\; = \;\; f(x, Z, \bar{Z}, \mu) &
  & \;\; = \;\; \omega(\mu) \; + \; \tilde f(x, Z, \bar{Z}, \mu), \\
  \dot{Z}_1 & \;\; = \;\; h_1(x, Z, \bar{Z}, \mu) &
  & \;\; = \;\; (\I \alpha_1(0) + \beta_1 + \I \delta_1) Z_1
  \; + \; \tilde h_1(x, Z, \bar{Z}, \mu), \\
  \dot{Z}_2 & \;\; = \;\; h_2(x, Z, \bar{Z}, \mu) &
  & \;\; = \;\; (\I \alpha_2(0) + \beta_2 + \I \delta_2) Z_2
  \; + \; \tilde h_2(x, Z, \bar{Z}, \mu)\, ,
\end{alignat*}
where $\mu_1 = \delta_1$, $\mu_2 = \delta_2$, $\mu_3 = \beta_1$
and $\mu_4 = \beta_2$ while the mute parameter~$\nu$ has been
dropped.

\begin{theorem}[Normal Form]
  \label{thm:normalform}
  Let $\ell \in \Z^2$ be such that $0 < \ell_1 < \ell_2$.
  Consider the real analytic vector field~\eqref{eq:X} where for
  $\kappa > n - 1$ and $\Gamma > 0$ the frequency vector
  $(\omega, \alpha)$ at $\mu = 0$ satisfies the Diophantine
  conditions
  \begin{equation}\label{eq:diophantine}
    \left|
      \inprod{k'}{\omega(0)} \; + \; \inprod{\ell'}{\alpha(0)}
    \right|
    \;\; \geq \;\; \frac{\Gamma}{(|k'| + |\ell'|)^\kappa}
  \end{equation}
  for all $(k', \ell') \in \Z^n \times \Z^2$ except for
  $(k', \ell') = (0, \ell^{\perp})$ and its integer multiples.
  For given order $M \in \N$ of normalisation there is a
  diffeomorphism
  \begin{align*}
    \Phi \: : \:\: \T^n \times \C^2 \times \R^s
    & \longrightarrow \T^n \times \C^2 \times \R^s,  
  \end{align*}
  real analytic in $x$, $Z$, $\bar{Z}$ and $\mu$, 
  close to the identity at $(Z, \mu) = (0,  0)$, such that the
  following holds.

  \samethought
  The vector field $X$ is transformed into normal form
  \begin{equation}\label{eq:normalform}
    \Phi_* X \;\; = \;\; N \; + \; R.
  \end{equation}
  The lower order part $N$ is
  \begin{align*}
    N & \;\; = \;\; \omega(0)\partial_x
        \; + \;
        (\I \alpha_1(0) + \beta_1 + \I \delta_1)
        Z_1 \partial_{Z_1}
        \; + \;
        (\I \alpha_2(0) + \beta_2 + \I \delta_2)
        Z_2 \partial_{Z_2}
        \\
      & \quad
        \; + \;
        \sum_m C^0_m
        \prod_{j=1}^4 \tau_j^{m_j}
        \prod_{k=1}^s \mu_{k}^{m_{4+k}}
        \partial_x \\
      & \quad
        \; + \;
        \sum_m C^1_m
        \prod_{j=1}^4 \tau_j^{m_j}
        \prod_{k=1}^s \mu_{k}^{m_{4+k}}
        Z_1\partial_{Z_1} \\
      & \quad
        \; + \;
        \sum_m C^2_m
        \prod_{j=1}^4 \tau_j^{m_j}
        \prod_{k=1}^s \mu_{k}^{m_{4+k}}
        Z_2\partial_{Z_2} \\
      & \quad
        \; + \; \text{c.c.},
  \end{align*}
  where c.c.\ is short for complex conjugate, 
  $C^j_m \in \C$ and $m \in \N^{4+s}$ with
  \[
    2 \;\; \le \;\; 2(m_1 + m_2) \; + \; |\ell| (m_3 + m_4)
    \; + \; \sum_{j=1}^s m_{4+j} \;\; \le \;\; M.
  \]
  The remainder $R$ satisfies the estimate
  \begin{align*}
    R \;\; = \;\;
      O_{M+1}(|Z|, |\mu|) \partial_x
    \; + \; O_{M+1}(|Z|, |\mu) Z_1 \partial_{Z_1}
    \; + \; O_{M+1}(|Z|, |\mu) Z_2 \partial_{Z_2}.
  \end{align*}
\end{theorem}

\newthought
The proof relies on a succession of normal form transformations,
see e.g.~\cite{BHW20}.
Taking $M = \max(4, |\ell|)$, the lower order part~$N$ of the
normal form can be written in the form 
\begin{align*}
  \dot{x} & \;\; = \;\; \omega(\mu)
             \; + \;
             \hat f(\tau_1, \tau_2,\mu)
             \; + \;
             a_0(\mu) \tau_3
             \; + \;
             c_0(\mu) \tau_4,
  \\
  \dot{Z}_j & \;\; = \;\; \left(
             \I \alpha_j(0) \, + \, \beta_j \, + \, \I \delta_j
             \, + \,
             p_j(\tau_1, \tau_2,\mu)
             \, + \,
             \I q_j(\tau_1, \tau_2,\mu)
             \right.
  \\
         & \qquad \left.
             \mbox{}+ \,
             [a_j(\mu) + \I b_j(\mu)] \tau_3
             \, + \,
             [c_j(\mu) + \I d_j(\mu)] \tau_4
             \right) Z_j, \quad j = 1, 2,
\end{align*}
where $\hat f(0, 0, \mu) = 0$ and
$p_j(0, 0, \mu) = q_j(0, 0, \mu) = 0$
while $\tau_3$ and~$\tau_4$ only
enter linearly --- that is what we
mean by lower order part. 

\samethought
The form of these equations already suggests that after discarding
the strong $1{:}1$ resonance there are still three distinct
situations.
\begin{itemize}
\item[-]
The first case is $|\ell| = 3$ --- the $1{:}2$~resonance ---
where the terms $\tau_3$ and~$\tau_4$ are larger than the terms
$\tau_1^2$ and~$\tau_2^2$.
\item[-]
The second case is $|\ell| = 4$ --- the $1{:}3$~resonance ---
where these terms are of the same order of magnitude.
\item[-]
The third case is $|\ell| \geq 5$ --- higher order resonances ---
where the terms $\tau_3$ and~$\tau_4$ are smaller than the terms
$\tau_1^2$ and~$\tau_2^2$.

\end{itemize}

\section{Analysis of the truncated normal form}
\label{ssec:weak:bifurcations}

In the normal form vector field~\eqref{eq:normalform} we truncate the
remainder $R$ so to obtain the $\T^n \times \T^1$ symmetric normal
form truncation~$N$.
Reducing the $\T^n$~symmetry we obtain a normal dynamics on~$\R^4$
that still has a $\T$~symmetry.
By abuse of notation we denote this reduced vector field by~$N$ as
well.

\samethought
Equilibria of the reduced~$N$ usually are called relative equilibria.
Observe that on $\T^n \times \R^4$ these correspond to invariant
$n$--tori.
In fact the normal system behaves like in the case where $n = 0$. 
For this reason the relative equilibria sometimes are referred to as
$0$--tori.

\samethought
Similarly periodic orbits of the reduced $N$ usually are called
relative periodic orbits.
Observe that on $\T^n \times \R^4$ these correspond to invariant
$(n+1)$--tori.
For this reason the relative periodic orbits sometimes are referred
to as $1$--tori.

\subsection{Passing to invariants}

Recall formula~\eqref{invariantphases} that relates $\tau_3$
and $\tau_4$ to $\tau_1, \tau_2$ and  the resonant angle
$\theta$.
Compared to polar co-ordinates, the invariants $\tau_1$
and~$\tau_2$ take the part of the radii, whereas $\tau_3$
and~$\tau_4$ take the part of $\cos \theta$ and~$\sin \theta$.

\samethought
Expressed in invariants, the normal form truncation~$N$
reads as 
\begin{align*}
  \dot{\tau}_1
  & \;\; = \;\; 2\left(
    \beta_1 + p_1 + a_1 \tau_3 + c_1 \tau_4
    \right) \tau_1
   , \\
  \dot{\tau}_2
  & \;\; = \;\; 2\left(
    \beta_2 + p_2 + a_2 \tau_3 + c_2 \tau_4
    \right) \tau_2, \\
  \dot{\tau}_3
  & \;\; = \;\; \left(
    \ell_2 (\beta_1 + p_1 + a_1 \tau_3 + c_1 \tau_4)
    \, + \,
    \ell_1 (\beta_2 + p_2 + a_2 \tau_3 + c_2 \tau_4)
    \right)\tau_3 \\
  & \quad
    - \left(
    \ell_2 (\delta_1 + q_1 + b_1 \tau_3 + d_1 \tau_4)
    \, - \,
    \ell_1 (\delta_2 + q_2 + b_2 \tau_3 + d_2 \tau_4)
    \right) \tau_4, \\
  \dot{\tau}_4
  & \;\; = \;\; \left(
    \ell_2 (\delta_1 + q_1 + b_1 \tau_3 + d_1 \tau_4)
    \, - \,
    \ell_1 (\delta_2 + q_2 + b_2 \tau_3 + d_2 \tau_4)
    \right)\tau_3 \\
  & \quad
    + \left(
    \ell_2 (\beta_1 + p_1 + a_1 \tau_3 + c_1 \tau_4)
    \, + \,
    \ell_1 (\beta_2 + p_2 + a_2 \tau_3 + c_2 \tau_4)
    \right) \tau_4.
\end{align*}
We introduce
\begin{equation*}
  A \; = \;
  \begin{pmatrix}
    a_1 & c_1 \\
    a_2 & c_2
  \end{pmatrix},
  \quad
  B \; = \;
  \begin{pmatrix}
    b_1 & d_1 \\ b_2 & d_2
  \end{pmatrix}.    
\end{equation*}
Preparing for a Taylor expansion of the right hand side, we
also introduce  
\begin{align*}
  p_{ij}(\mu)
  & \; = \; \dd {p_i} {\tau_j}(0, 0, \mu), \quad
  P(\mu) \; = \;
  \begin{pmatrix}
    p_{11}(\mu) & p_{12}(\mu) \\
    p_{21}(\mu) & p_{22}(\mu)
  \end{pmatrix}
  \\
  \intertext{as well as}
  q_{ij}(\mu)
  & \; = \; \dd {q_i} {\tau_j}(0, 0, \mu), \quad
  Q(\mu) \; = \;
  \begin{pmatrix}
    q_{11}(\mu) & q_{12}(\mu) \\
    q_{21}(\mu) & q_{22}(\mu)
  \end{pmatrix}.
\end{align*}
We zoom in on the origin of $\R^4$ by scaling time, the
invariants and the parameters as
\begin{align*}
  t & \; = \; \eps^{-2}\tilde t, \quad
  \tau_j \; = \; \eps^2 \sigma_j, \quad 
  \tau_{j+2} \; = \; \eps^{|\ell|} \psi_j, \quad
  \beta_j \;  = \; \eps^2 \gamma_j, \quad
  \delta_j \; = \; \eps^2 \eta_j, \quad j = 1, 2,
\end{align*}
thereby splitting $\tau$ into the amplitude~$\sigma$ and
the phase~$\psi$ (recall from~\eqref{invariantphases} that
$\psi$ is a scaled version of $\cos \theta$ and~$\sin \theta$).
This allows to introduce $\tilde{p}$ and~$\tilde{q}$ by
\begin{align*}
  p(\eps^2 \sigma, \mu)
  & \;\; = \;\; \eps^2 P(\mu) \sigma
    \; + \; \eps^4 \tilde{p}(\sigma, \eps^2, \mu), \\[3pt]
  q(\eps^2 \sigma, \mu)
  & \;\; = \;\; \eps^2 Q(\mu) \sigma
    \; + \; \eps^4 \tilde{q}(\sigma, \eps^2, \mu).
\end{align*}
Finally, the detuning $\delta$ of the frequencies now results in
the parameter
\begin{equation}\label{eq:frequency}
  \xi
  \;\; = \;\; \inprod{\ell^{\perp}}{\eta}
  \;\; = \;\; \eps^{-2} \inprod{\ell^{\perp}}{\delta}
\end{equation}
that detunes the frequency ratio.

\samethought
Introduce $\ell^* := (\ell_2, \ell_1)$, and, for a vector
$y \in \R^{\nu}$, the $\nu \times \nu$ matrix $\diag(y)$ as
the diagonal matrix whose $i$'th diagonal element is $y_i$. 
The equations of motion defined by $N$ can then be written
as 
\begin{subequations}
  \label{eq:sigma-psi}
  \begin{align}
    \label{eq:sigma-psi:sigma}
    \dot{\sigma}
    & \;\; = \;\; 2 \diag(\gamma + P\sigma
           + \eps^2 \tilde p + \eps^{|\ell|-2}A \psi)
      \sigma, \\ 
    \label{eq:sigma-psi:psi}
    \dot{\psi}
    & \;\; = \;\;
      \begin{pmatrix}
        \inprod{\ell^*}{\gamma + P \sigma + \eps^2 \tilde p +
          \eps^{|\ell|-2}A \psi}
        & -\xi - \inprod{\ell^{\perp}}{Q
          \sigma + \eps^2 \tilde q + \eps^{|\ell|-2}B \psi}
        \\
        \xi + \inprod{\ell^{\perp}}{Q \sigma + \eps^2 \tilde q +
          \eps^{|\ell|-2}B \psi} & 
        \inprod{\ell^*}{\gamma + P \sigma + \eps^2 \tilde p +
          \eps^{|\ell|-2}A \psi}
      \end{pmatrix}
      \psi .
  \end{align}
\end{subequations}

\samethought
The intricacies of the resonance at hand are encoded
in~\eqref{eq:sigma-psi:psi}.
The syzygy takes the form
\begin{equation}
  \label{eq:sigma-psi:syzygy}
  \sigma_1^{\ell_2} \sigma_2^{\ell_1} \; - \;
  G_\ell (\psi_1^2 + \psi_2^2) \;\; = \;\; 0.
\end{equation}
The equations of motions defined by~$N$ are then defined on the
reduced phase space
\begin{equation*}
  \sP \;\; = \;\;
  \left\{ (\sigma,\psi) \in \R^2 \times \R^2 \,:\,
    \sigma_1\ge 0,\; \sigma_2 \ge 0,\; 
    \psi_1^2 + \psi_2^2 \; = \;
    \frac{\sigma_1^{\ell_2} \sigma_2^{\ell_1}}{G_\ell}
  \right\},
\end{equation*}
which is a semi-algebraic variety.
Geometrically, $\sP$ is a degenerate circle bundle over the
basis $\sB = \R_{\ge 0}^2$ that is degenerate at the boundary
of~$\sB$.  

\samethought
As we are not trying to perform an exhaustive analysis, we
restrict the analysis to the situation where there is a
subordinate Hopf bifurcation from a $2$--dimensional to a
$3$--dimensional torus. 
In particular, we investigate the equations of
motion~\eqref{eq:sigma-psi} under the assumptions
that $p_{11} > 0$, $p_{22} < 0$ and $\det P > 0$.
To fix thoughts we restrict to $n = 0$.

\subsection{Equilibria on the basis~$\sB$}

For $\eps = 0$ the system \eqref{eq:sigma-psi} simplifies to
\begin{subequations}
  \label{eq:sigma-psi-approx}
  \begin{align}
    \label{eq:sigma-psi-approx:sigma}
    \dot{\sigma} & \;\; = \;\; 2 \diag(\gamma + P \sigma)
                   \sigma, \\
    \label{eq:sigma-psi-approx:psi}
    \dot{\psi} & \;\; = \;\; \begin{pmatrix}
      \inprod{\ell^*}{\gamma + P \sigma}
      &
      - \xi - \inprod{\ell^{\perp}}{Q \sigma}
      \\
      \xi + \inprod{\ell^{\perp}}{Q \sigma}
      &
      \inprod{\ell^*}{\gamma + P \sigma}
    \end{pmatrix}
    \psi.
  \end{align}
\end{subequations}
Remark that the system~\eqref{eq:sigma-psi-approx} is skew:
the basis dynamics~\eqref{eq:sigma-psi-approx:sigma}, which
is defined on the basis $\sB$ of $\sP$, is decoupled from the
fibre dynamics \eqref{eq:sigma-psi-approx:psi} and, in fact,
drives the fibre dynamics.
Also note that the basis dynamics has Lotka--Volterra structure;
that is, the components of the vector field are quadratic
polynomials and the axes $\sigma_1 = 0$ and $\sigma_2 = 0$ are
invariant. 
Finally, note that~\eqref{eq:sigma-psi:sigma} is the same
as in the non-resonant case, compare with~\cite{GH83, kuz95, li16}.

\newthought
Equilibria of the basis dynamics are the central equilibrium 
\begin{equation*}
  \bar{\sigma} \;\; = \;\; 0,
\end{equation*}
which corresponds to the equilibrium of the original system
on $\R^4$, the boundary equilibria
\begin{equation*}
  \bar{\sigma} \; = \;
  \begin{pmatrix}
    -\gamma_1/p_{11} \\
    0
  \end{pmatrix},\quad 
  \bar{\sigma} \; = \;
  \begin{pmatrix}
    0 \\
    -\gamma_2/p_{22}
  \end{pmatrix}, \\
\end{equation*}
which correspond to periodic orbits in the original system,
and the interior equilibrium
\begin{equation*}
  \bar{\sigma} \;\; = \;\; - P^{-1} \gamma.
\end{equation*}
which corresponds to an invariant $2$--torus.
All equilibria are defined for those values of $\gamma$ such
that $\bar \sigma_1 \geq 0$ and $\bar{\sigma}_2 \geq 0$.

\subsection{Hopf bifurcations on the basis~$\sB$}

Linearising the flow to investigate the stability of the
equilibria yields
\begin{equation*}
  \dot{\sigma} \;\; = \;\;
  \begin{pmatrix}
    \gamma_1 + 2 p_{11} \bar{\sigma}_1 + p_{12} \bar{\sigma}_2
    & p_{12} \bar{\sigma}_1 \\
    p_{21} \bar \sigma_2
    & \gamma_2 + p_{21} \bar{\sigma}_1 + 2 p_{22} \bar{\sigma}_2
  \end{pmatrix}
  (\sigma - \bar{\sigma}) \; + \; O(\|\sigma - \bar{\sigma}\|^2).
\end{equation*}
For the central equilibrium, this reduces to
\begin{equation*}
  \dot{\sigma} \;\; = \;\;
  \begin{pmatrix}
    \gamma_1 & 0 \\
    0 & \gamma_2 
  \end{pmatrix}
  \sigma \; + \; O(\|\sigma\|^2).
\end{equation*}
The central equilibrium is attracting if $\gamma_j < 0$ for
$j = 1, 2$ and repelling if $\gamma_j > 0$ for $j = 1, 2$.
Correspondingly, the curves
\begin{equation*}
  H^a_{0 \to 1} \quad : \quad \gamma_1 \; = \; 0
  \qquad \mbox{and} \qquad
  H^b_{0 \to 1} \quad : \quad \gamma_2 \; = \; 0
\end{equation*}
are Hopf bifurcation curves for the original system, where a
$0$--torus changes stability.

\begin{figure}[ht]
  \centering
  \begin{minipage}[b]{0.24\linewidth}
    \centering
    \begin{minipage}[b]{0.7\linewidth}
      \centering
      \includegraphics[width=\textwidth]{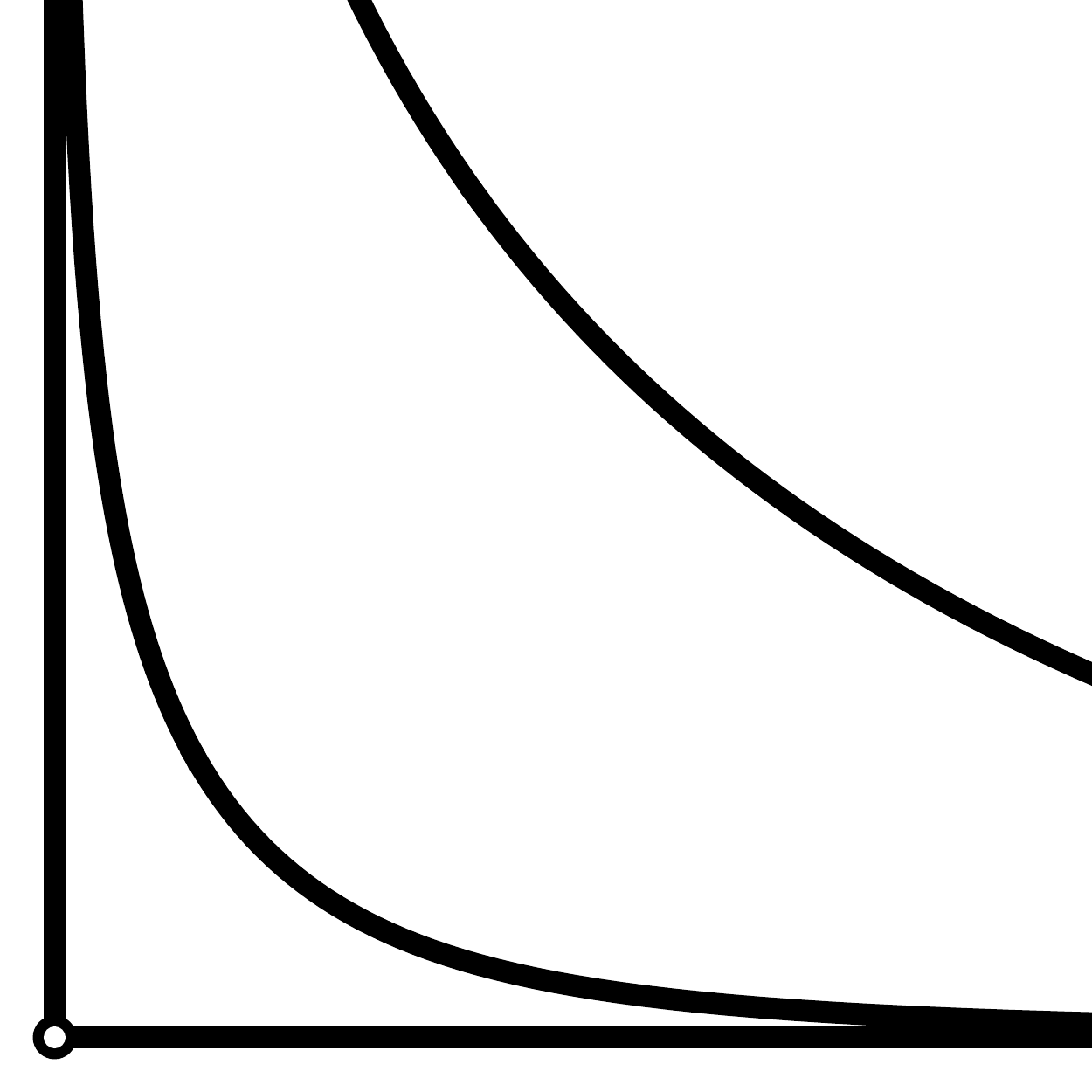}
      \subcaption{Region I}
    \end{minipage}
    \begin{minipage}[b]{0.7\linewidth}
      \includegraphics[width=\textwidth]{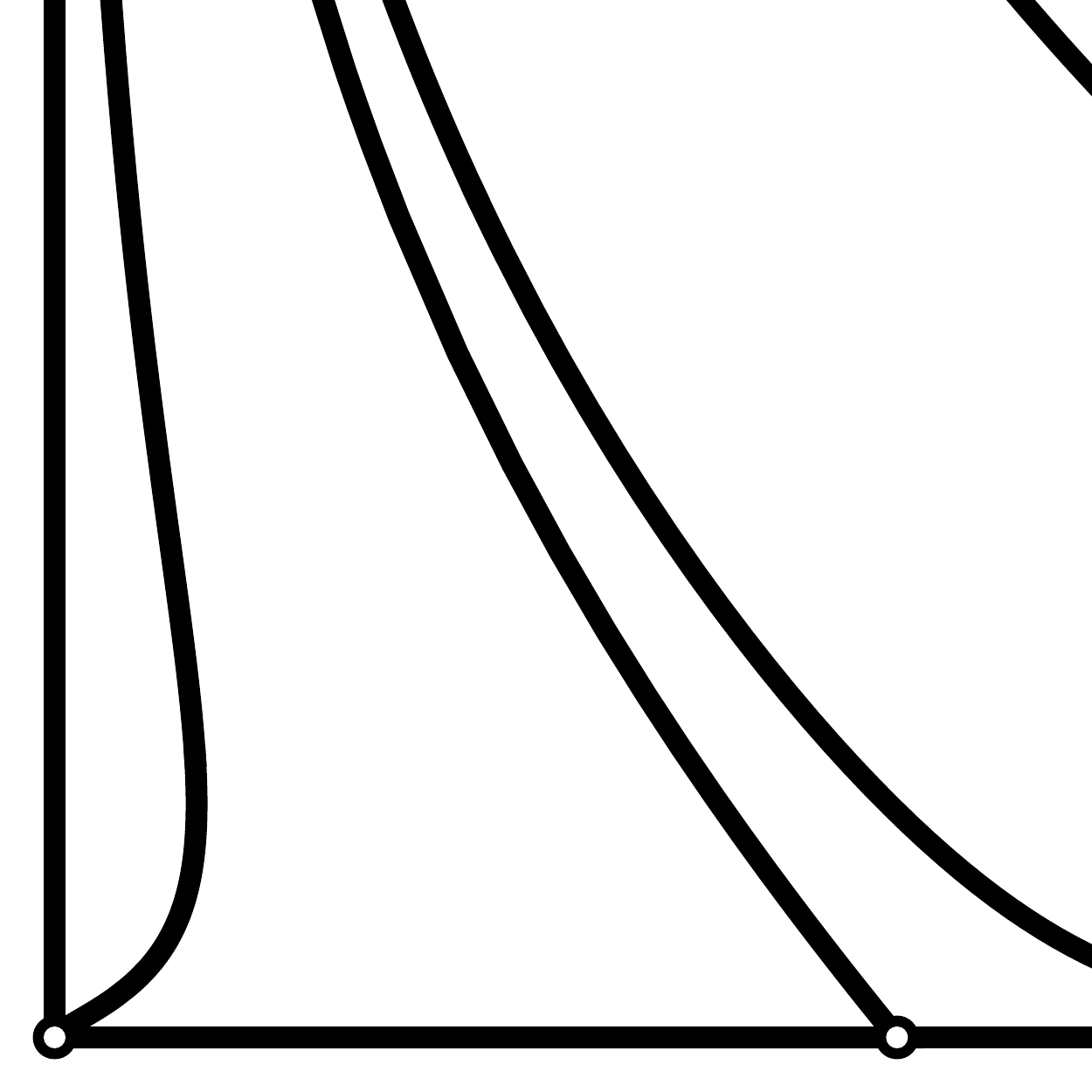}
      \subcaption{Region II}
    \end{minipage}
  \end{minipage}
  \hfill
  \begin{minipage}[b]{0.5\linewidth}
    \includegraphics[width=\textwidth]{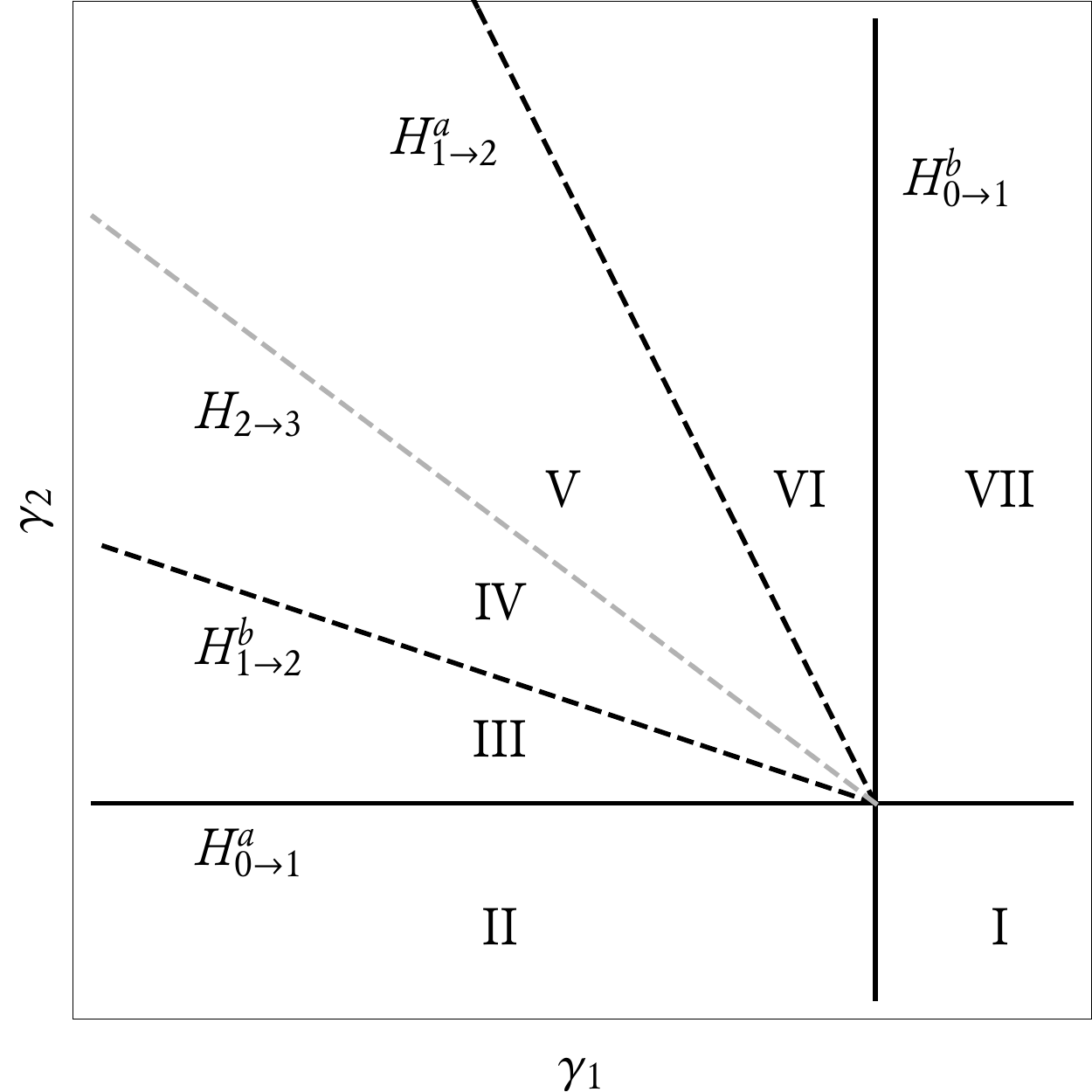}  
  \end{minipage}
  \hfill
  \begin{minipage}[b]{0.24\linewidth}
    \centering
    \begin{minipage}[b]{0.7\linewidth}
      \centering
      \includegraphics[width=\textwidth]{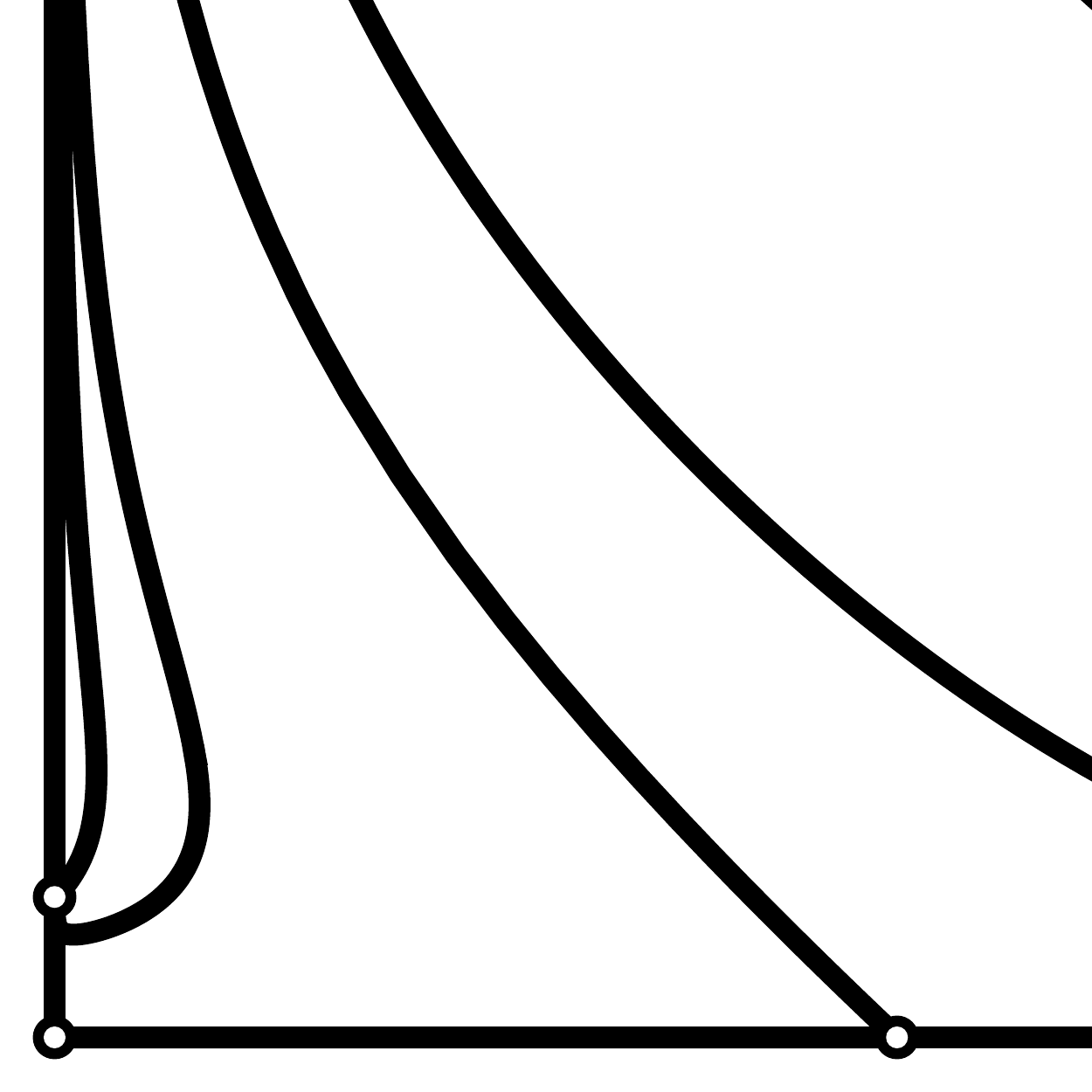}
      \subcaption{Region III}
    \end{minipage}
    \begin{minipage}[b]{0.7\linewidth}
      \includegraphics[width=\textwidth]{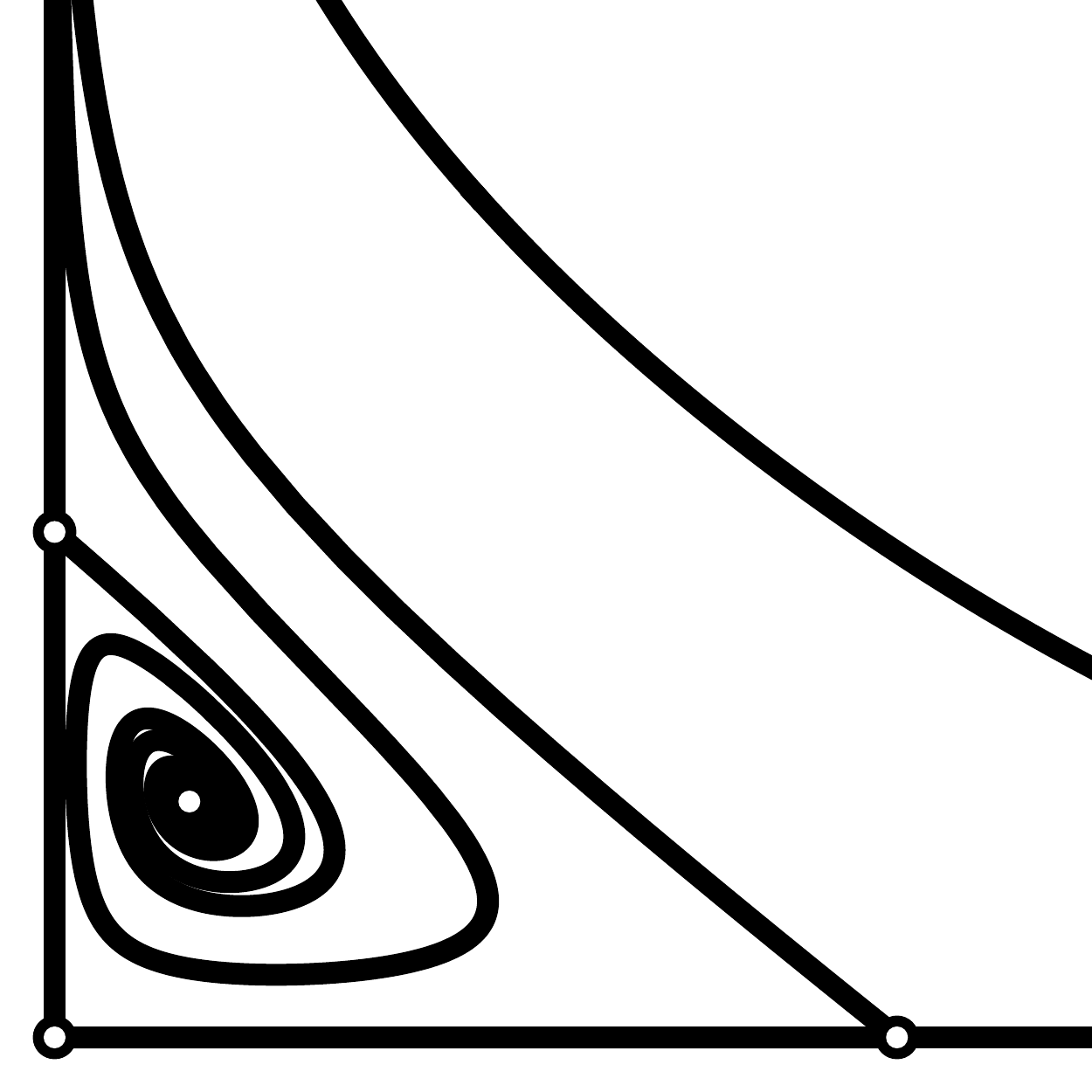}
      \subcaption{Region IV}
    \end{minipage}
  \end{minipage}
  \begin{minipage}[b]{\linewidth}
    \centering
    \begin{minipage}[b]{0.168\linewidth}
      \centering
      \includegraphics[width=\textwidth]{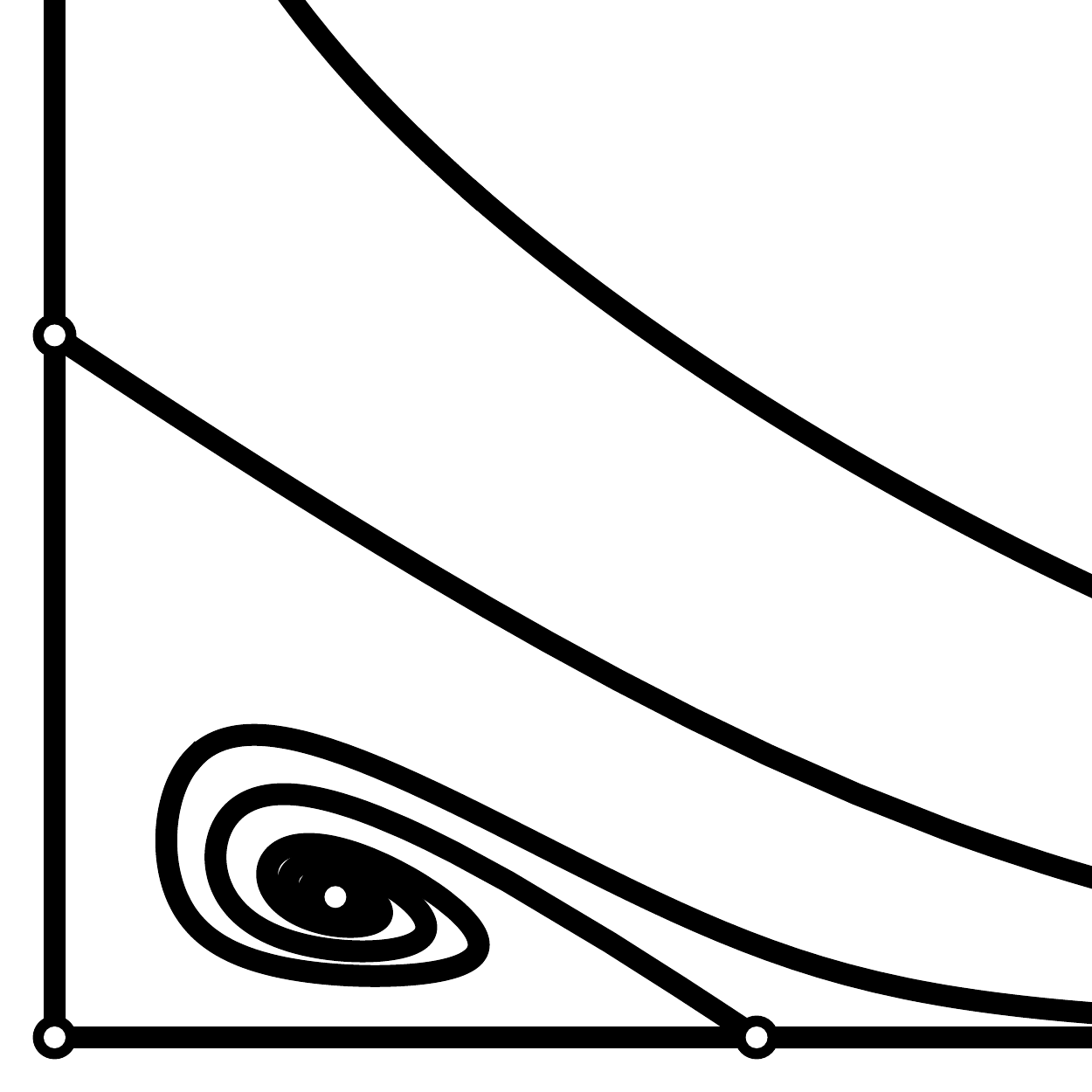}
      \subcaption{Region V}
    \end{minipage}
    \hspace{0.05\textwidth}
    \begin{minipage}[b]{0.168\linewidth}
      \centering
      \includegraphics[width=\textwidth]{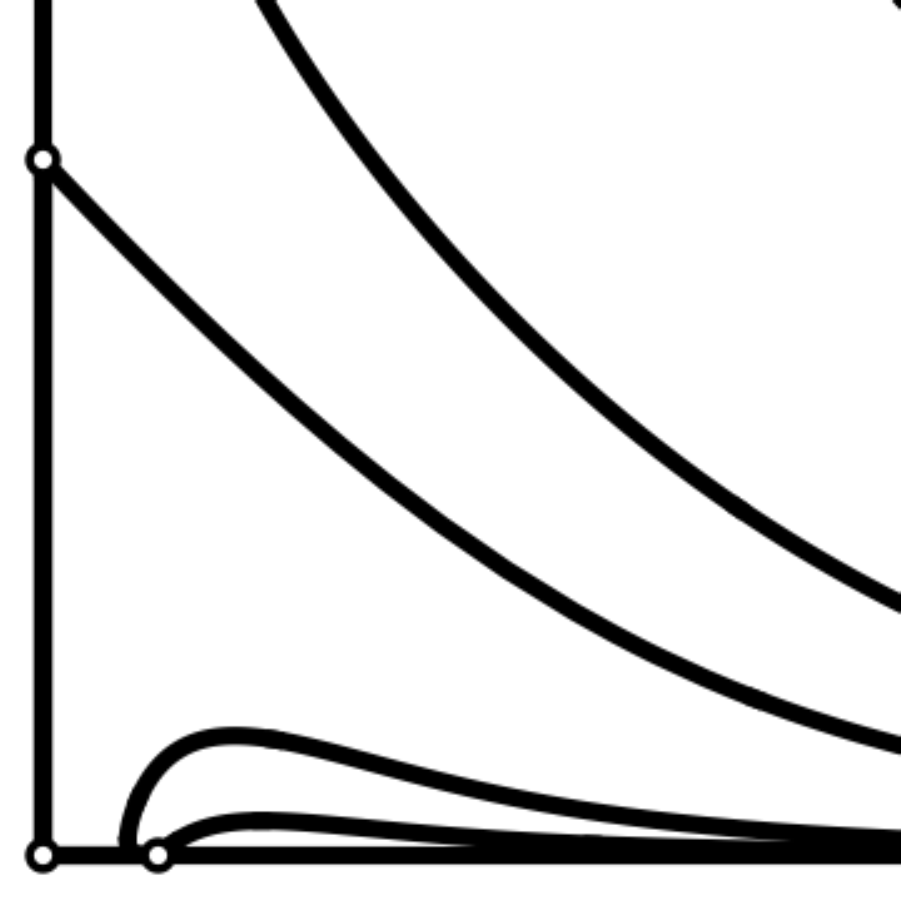}
      \subcaption{Region VI}
    \end{minipage}
    \hspace{0.05\textwidth}
    \begin{minipage}[b]{0.168\linewidth}
      \centering
      \includegraphics[width=\textwidth]{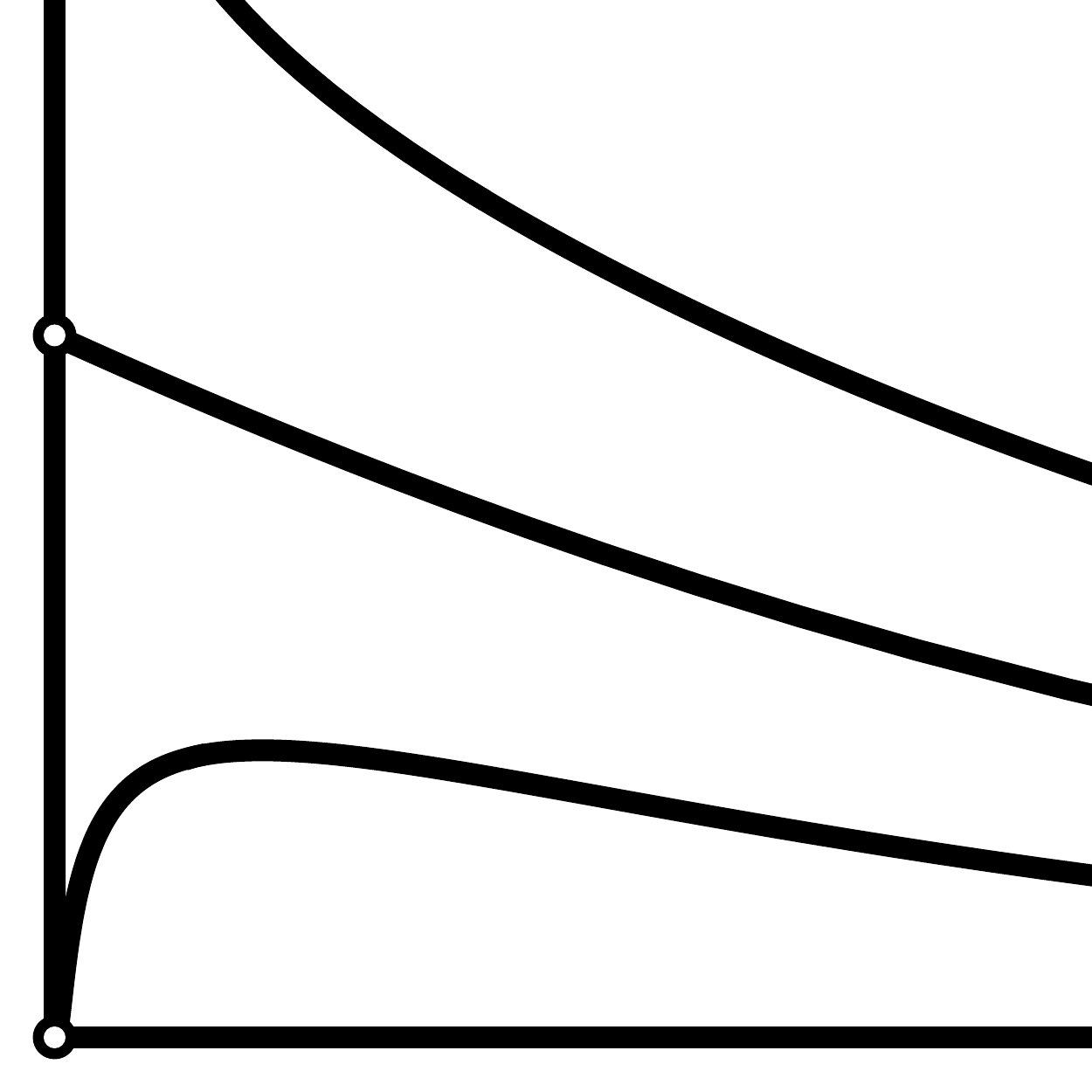}
      \subcaption{Region VII}
    \end{minipage}
    \hspace{0.05\textwidth}
    \begin{minipage}[b]{0.168\linewidth}
      \centering
      \includegraphics[width=\textwidth]{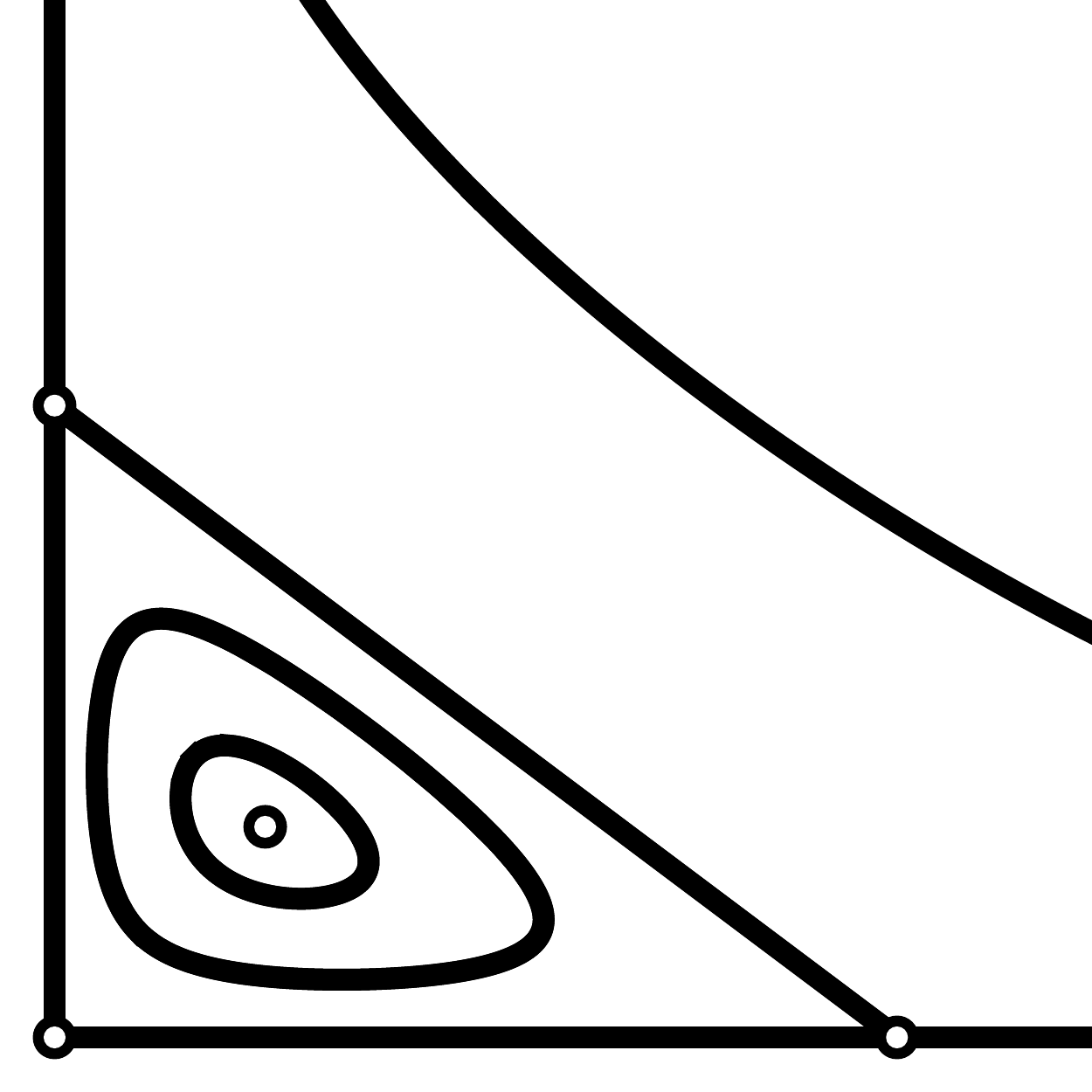}
      \subcaption{$H_{2\to 3}$}
    \end{minipage}
  \end{minipage}
  \caption{Bifurcation diagram of~\eqref{eq:sigma-psi-approx} 
    with $P = (^{2^{-1}}_{-1} {}^{\;\; 3}_{-1})$.
  The parameters $\gamma_1$ and $\gamma_2$ are scaled versions of the
  real parts $\beta_j$ of the two pairs of eigenvalues.  
  The regions I--VII in the central picture refer to the surrounding phase
  portraits on the basis~$\sB$.
  For the coding of the subordinate Hopf bifurcations see the main text.}
  \label{fig:hopf}  
\end{figure}

\samethought
The $1$--tori in $\R^4$ bifurcating off from the origin correspond
to the boundary equilibria.
These tori pass in turn through a secondary Hopf bifurcation.
For the boundary equilibrium
$\bar{\sigma} = (-\gamma_1 / p_{11}, 0)$,
which under our assumptions on~$P$ exists if $\gamma_1 > 0$, we
obtain    
\begin{equation*}
  \dot{\sigma} \;\; = \;\;
  \begin{pmatrix}
    -\gamma_1 & - \frac{p_{21}}{p_{11}} \gamma_1 \\
    0 & \gamma_2 - \frac{p_{21}}{p_{11}} \gamma_1
  \end{pmatrix}
  (\sigma - \bar \sigma) \; + \; O(\|\sigma - \bar \sigma\|^2).
\end{equation*}
Hence, this equilibrium changes from attractor to saddle at
\begin{equation*}
  H^a_{1 \to 2} \quad : \quad p_{21} \gamma_1 \; - \; p_{11} \gamma_2
  \;\; = \;\; 0.
\end{equation*}
Similarly, the boundary equilibrium
$\bar{\sigma} = (0, -\gamma_2 / p_{22})$ changes from attractor to
saddle at 
\begin{equation*}
  H^b_{1 \to 2} \quad : \quad p_{22} \gamma_1 \; - \; p_{12} \gamma_2
  \;\; = \;\; 0.
\end{equation*}
Linearising the flow at the interior equilibrium yields
\begin{equation*}
  \dot{\sigma} \;\; = \;\;
  \begin{pmatrix}
    p_{11} \bar{\sigma}_1 & p_{12} \bar{\sigma}_1 \\
    p_{21} \bar{\sigma}_2 & p_{22} \bar{\sigma}_2
  \end{pmatrix}
  (\sigma - \bar{\sigma}) \; + \; O(\|\sigma - \bar{\sigma}\|^2).
\end{equation*}
The characteristic equation of the matrix on the right hand
side is
\begin{align*}
  0 & \;\; = \;\; \det
      \begin{pmatrix}
        p_{11} \bar{\sigma}_1 - \lambda & p_{12} \bar{\sigma}_1 \\
        p_{21} \bar{\sigma}_2 & p_{22} \bar{\sigma}_2 - \lambda 
      \end{pmatrix}
      \;\; = \;\; \bar{\sigma}_1 \bar{\sigma}_2
        \det
      \begin{pmatrix}
        p_{11} - \lambda/\bar{\sigma}_1 & p_{12}  \\
        p_{21} \bar{\sigma}_2 & p_{22} -\lambda/\bar{\sigma}_2 
      \end{pmatrix}
  \enspace .
\end{align*}
The equilibrium has imaginary eigenvalues if
\begin{equation*}
  \text{trace}
  \begin{pmatrix}
    p_{11} \bar{\sigma}_1 & p_{12} \bar{\sigma}_1 \\
    p_{21} \bar{\sigma}_2 & p_{22} \bar{\sigma}_2
  \end{pmatrix}
  \;\; = \;\;
  p_{11} \bar{\sigma}_1 + p_{22} \bar{\sigma}_2 
  \;\; = \;\; 0 
\end{equation*}
and
\begin{equation*}
  \det \begin{pmatrix}
    p_{11} \bar{\sigma}_1 & p_{12} \bar{\sigma}_1 \\
    p_{21} \bar{\sigma}_2 & p_{22} \bar{\sigma}_2
  \end{pmatrix}
  \;\; = \;\;
  (\det P) \bar{\sigma}_1 \bar{\sigma}_2 \;\; > \;\; 0.
\end{equation*}
The first condition can only be satisfied if the signs of
$p_{11}$ and~$p_{22}$ are opposite, as we have assumed.
The second condition is always satisfied given our assumption
that $\det P > 0$.
The first condition can be written as
\begin{equation*}
  H_{2 \to 3}
  \quad : \quad
  p_{22} (p_{21} - p_{11}) \gamma_1
  \; + \;
  p_{11} (p_{12} - p_{22}) \gamma_2 \;\; = \;\; 0.
\end{equation*}
It is well-known~\cite{GH83} that this Hopf bifurcation is
degenerate, and that the base system has a first integral at
the $H_{2 \to 3}$ bifurcation value.
Using the condition, the normal frequency at Hopf bifurcation can
be expressed as a function of the scaled parameter~$\gamma_1$ as
\begin{equation*}
  \sqrt{(\det P) \bar{\sigma}_1 \bar{\sigma}_2}
  \;\; = \;\;
  |\gamma_1| |p_{21} - p_{11}|
  \sqrt{(p_{12}p_{21} - p_{11} p_{22})
    p_{11}/p_{22}} 
\end{equation*}
The location of the Hopf bifurcation curves is illustrated in
figure~\ref{fig:hopf}.

\subsection{Dynamics on the reduced phase space~$\sP$}

Joining the fibre dynamics \eqref{eq:sigma-psi-approx:psi} to
the basis dynamics~\eqref{eq:sigma-psi-approx:sigma}, for
$\eps=0$ we find two degeneracies.
First, for parameters at the $H_{2\to 3}$ Hopf bifurcation
value, the basis dynamics has a first integral.
Second, an equilibrium~$\bar{\sigma}$ of the basis dynamics
corresponds to a limit cycle as long as the frequency
$\xi + \inprod{\ell^{\perp}}{Q\bar{\sigma}}$ does not vanish and
to a circle of equilibria where it does vanish.
In the latter case we speak of a resonance droplet. 

\begin{figure}[ht]
\begin{center}
\begin{picture}(450,200)
   \put(0,0){\includegraphics[width=200pt,keepaspectratio]{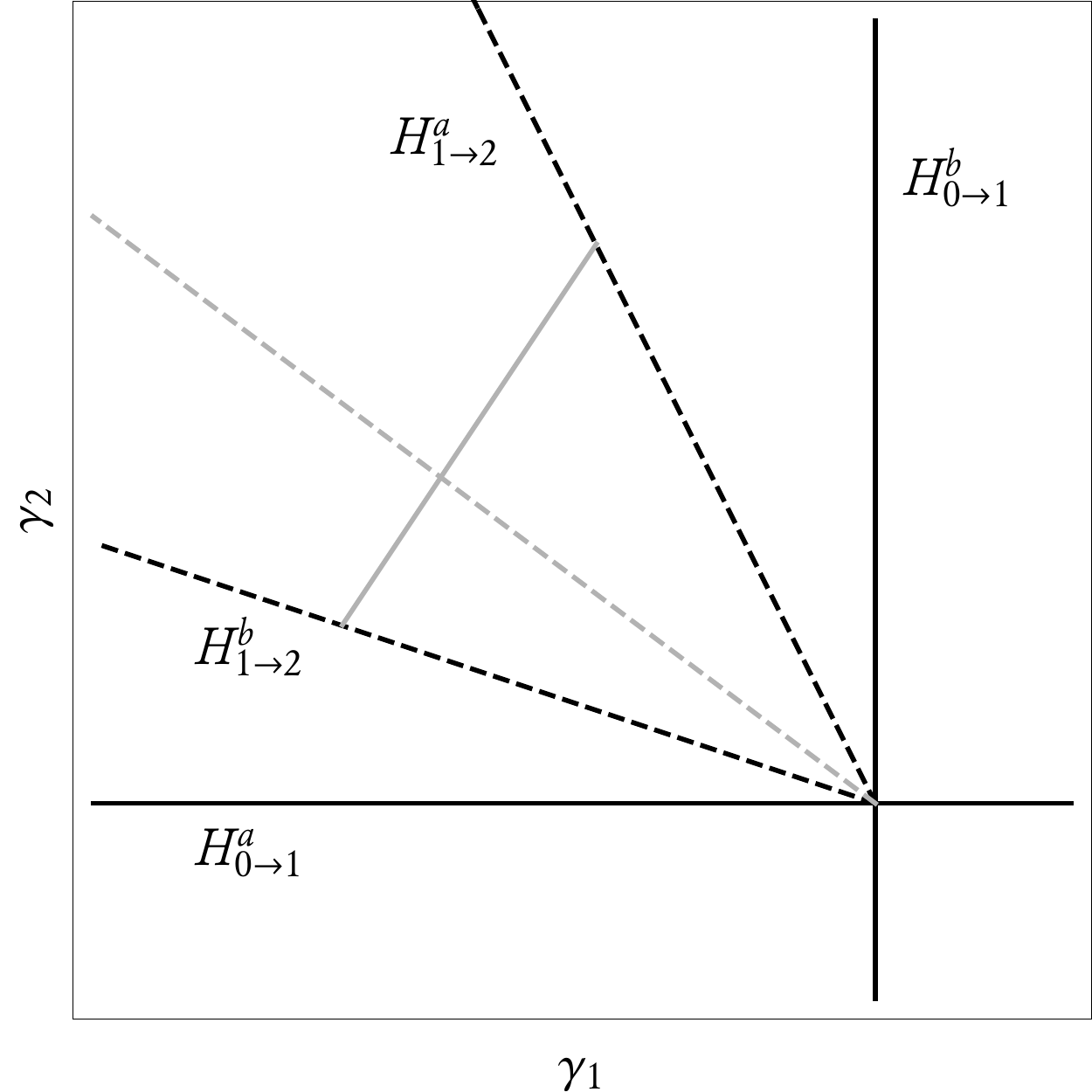}}
   \put(250,30){\includegraphics[width=200pt,keepaspectratio]{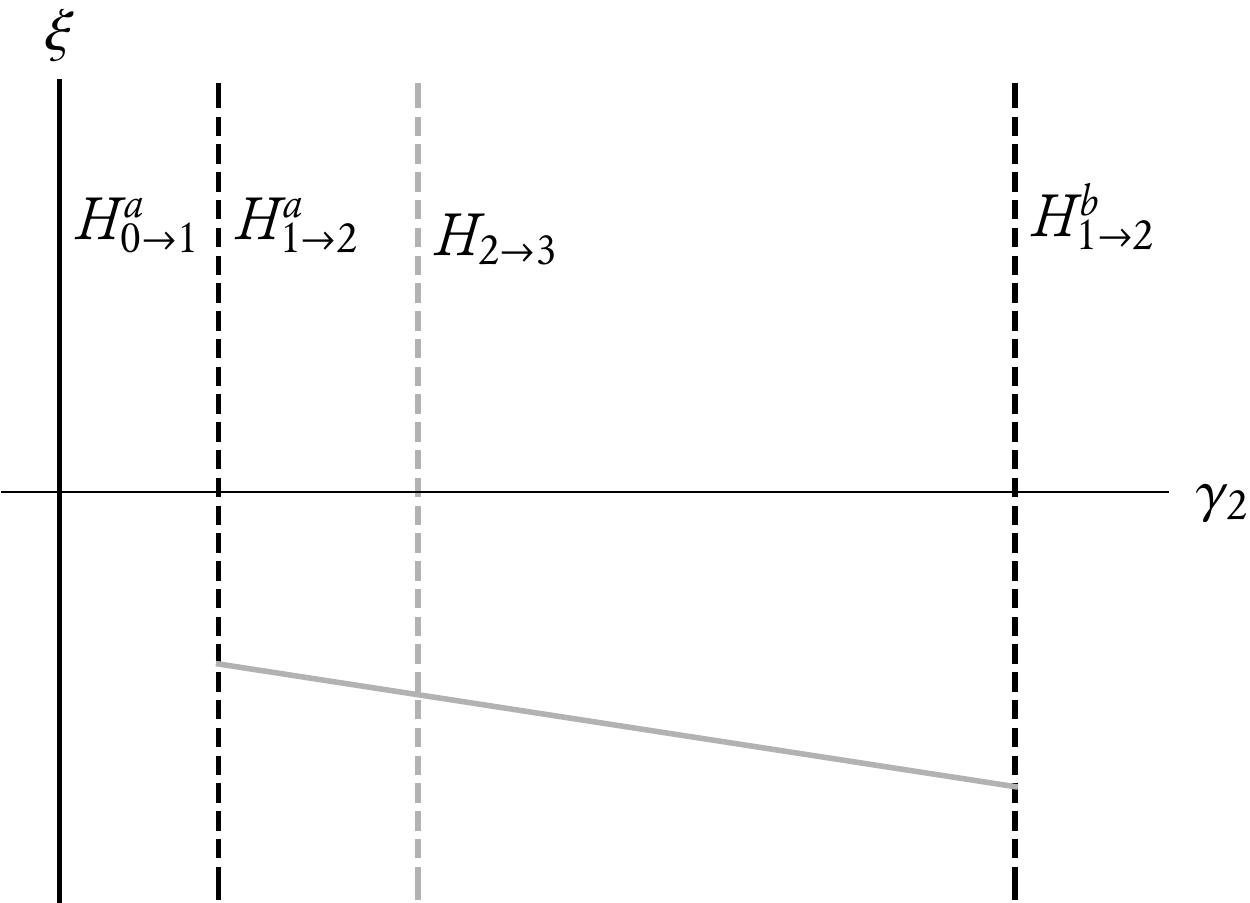}}
\end{picture}
\end{center}
\caption{The normal form~\eqref{eq:sigma-psi-approx} does
  not fully resolve the resonance droplet, which is shown as a grey segment.
  In the left hand picture we use the scaled 
    eigenvalue parameters $\gamma_1$ and~$\gamma_2$. 
    In the right hand picture we replace $\gamma_1$
    by the frequency~$\xi,$ defined in \eqref{eq:frequency}.
  Also not resolved is the Hopf line~$H_{2 \to 3}$, which
  should yield a whole range of parameters with periodic
  orbits on the basis~$\sB$, which is also shown in grey.
}\label{fig:degenerate}
\end{figure}

\newthought
Now consider~\eqref{eq:sigma-psi} for $\eps > 0$.
There are two situations, according to whether the frequency
$\xi + \inprod{\ell^{\perp}}{Q \bar{\sigma}}$ of the fibre dynamics
does or does not vanish.
If $\bar\sigma$ is a hyperbolic equilibrium of the basis
dynamics and if $\xi + \inprod{\ell^{\perp}}{Q \bar{\sigma}}$ does
not vanish, the limit cycle $\{ \bar{\sigma} \} \times \T$ survives a
small perturbation.
If $\bar\sigma$ is at the subordinate Hopf bifurcation~$H_{2 \to 3}$
we have to normalise to higher order than~$4$ to resolve the
degeneracy.
As shown in~\cite{GH83} already one additional order suffices.
This has consequences only for the low order $1{:}2$
and $1{:}3$~resonances where indeed $|\ell| \leq 4$.
For the higher order resonances we already do normalise to order
$|\ell| \geq 5$.
Note that next to the Hopf line~$H_{2 \to 3}$ also the resonance
bubbles are not fully resolved, see figure~\ref{fig:degenerate}.

\samethought
In this way the family of periodic orbits at a single parameter value
gets re-distributed over a range of parameters, starting at the Hopf
bifurcation~$H_{2 \to 3}$ and ending at the parameter value~$\Het$
of a heteroclinic bifurcation where a closed orbit disappears
in a blue sky bifurcation.

\subsection{Saddle-node bifurcations on the reduced phase space~$\sP$}

\samethought
We proceed to find the equilibria of~\eqref{eq:sigma-psi}, which
satisfy the equations
\begin{subequations}
\label{eq:sigma-psi:isoclines}
  \begin{align}
    0 & \;\; = \;\; \gamma \; + \; P \sigma
        \; + \; \eps^2 \tilde{p} \; + \; \eps^{|\ell|-2} A \psi,
  \label{eq:sigma-psi:isoclines:a} \\
    0 & \;\; = \;\; \xi \; + \; \inprod{\ell^{\perp}}{Q \sigma
        \; + \; \eps^2 \tilde{q} \; + \; \eps^{|\ell|-2} B \psi},
  \label{eq:sigma-psi:isoclines:b}
  \end{align}
\end{subequations}
together with the syzygy~\eqref{eq:sigma-psi:syzygy}.

\begin{figure}[htp]
  \centering
  \includegraphics[height=0.5\textheight]{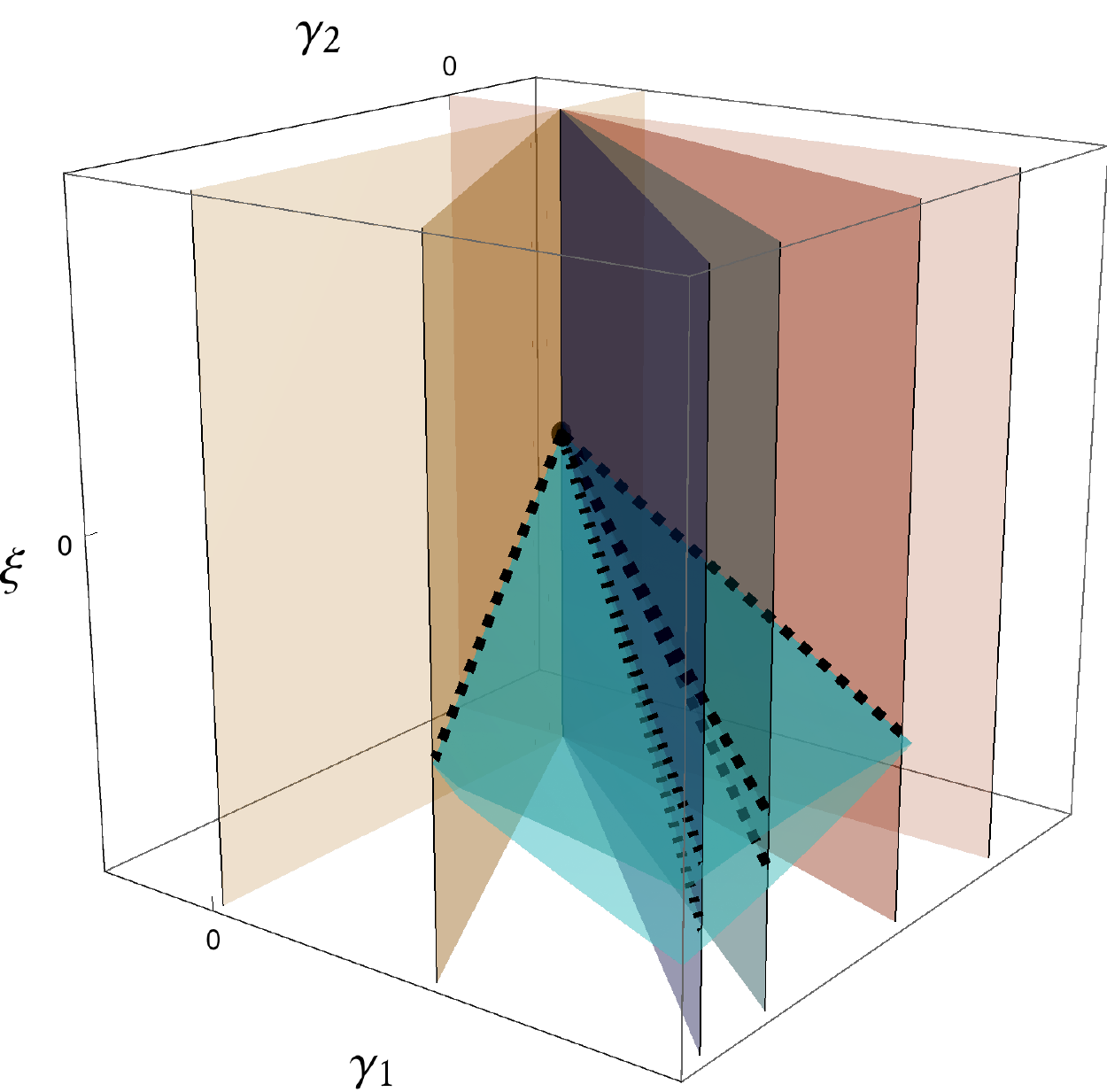}
  \caption{Local bifurcation diagram of the family~\eqref{eq:sigma-psi}
    in the space of the scaled eigenvalue parameters
    $\gamma_1$, $\gamma_2$ and the frequency~$\xi,$ defined
    in \eqref{eq:frequency}. 
    The orange planes are the Hopf bifurcations taken over from
    figure~\ref{fig:hopf}. 
    In green we furthermore include the saddle-node bifurcations
    originating from the resonance droplets in
    figure~\ref{fig:degenerate}.
    Finally, the purple plane denotes the heteroclinic
    bifurcations, which originate from resolving the
    $H_{2\to 3}$ bifurcation.
  } 
  \label{fig:bd}
\end{figure}

\newthought
As $P$ is invertible, which is the only situation we consider,
equation~\eqref{eq:sigma-psi:isoclines:a} can be solved in the form
\begin{equation*}
  \sigma
  \;\; = \;\;\Phi_0(\gamma,\eps^2)
  \; + \;
  \eps^{|\ell|-2}
  \Phi_1(\gamma, \eps^2, \eps^{|\ell|-2} A \psi) A \psi
\end{equation*}
where $\Phi_0(\gamma, 0) = -P^{-1} \gamma$ and
$\Phi_1(\gamma, 0, 0)= -P^{-1}$. 
Substitution into~\eqref{eq:sigma-psi:isoclines:b} yields
\begin{align}
\label{eq:sigma-psi:isocline:b:variant}
  0 & \;\; = \;\; \xi
      \; + \; \inprod{\ell^{\perp}}{Q\Phi_0 + \eps^2 \tilde q}
      \; + \; \eps^{|\ell|-2}
        \inprod{\ell^{\perp}}{(B+\Phi_1 A)\psi}.
\end{align}
To describe the saddle-node bifurcation curves, we perform a
blow-up by introducing the scaled parameter $\zeta$ as
\begin{equation*}
  \eps^{ |\ell|-2} \zeta
  \;\; = \;\;
  \xi
  \; + \;
  \inprod{\ell^{\perp}}{Q\Phi_0 + \eps^2 \tilde q}.
\end{equation*}
Dividing by~$\eps^{|\ell|-2}$, this
transforms~\eqref{eq:sigma-psi:isocline:b:variant} into
\begin{align}
  \label{eq:sigma-psi:resonance}
  0 & \;\; = \;\; \zeta \; + \;
  \inprod{\ell^{\perp}}{(B+\Phi_1 A)\psi}.
\end{align}
As an equation in~$\psi$, \eqref{eq:sigma-psi:resonance}~describes
a family of curves that for $\eps = 0$ reduces to a family of lines,
since $\Phi_1$, $A$ and~$B$ are independent of~$\psi$.
For fixed values of $\sigma$, the syzygy~\eqref{eq:sigma-psi:syzygy}
defines a circle.
The curves~\eqref{eq:sigma-psi:resonance} intersect the circle
in either zero or two points, or they touch the circle.
In the latter situation, we have a saddle-node bifurcation of
equilibria, defining the boundaries of a resonance droplet of
width $\eps^{|\ell| - 2}$.
The resulting bifurcation diagram is given in figure~\ref{fig:bd}.

\begin{figure}[ht]
  \centering
  \begin{minipage}[b]{0.425\linewidth}
    \centering
    \includegraphics[width=\textwidth]{%
      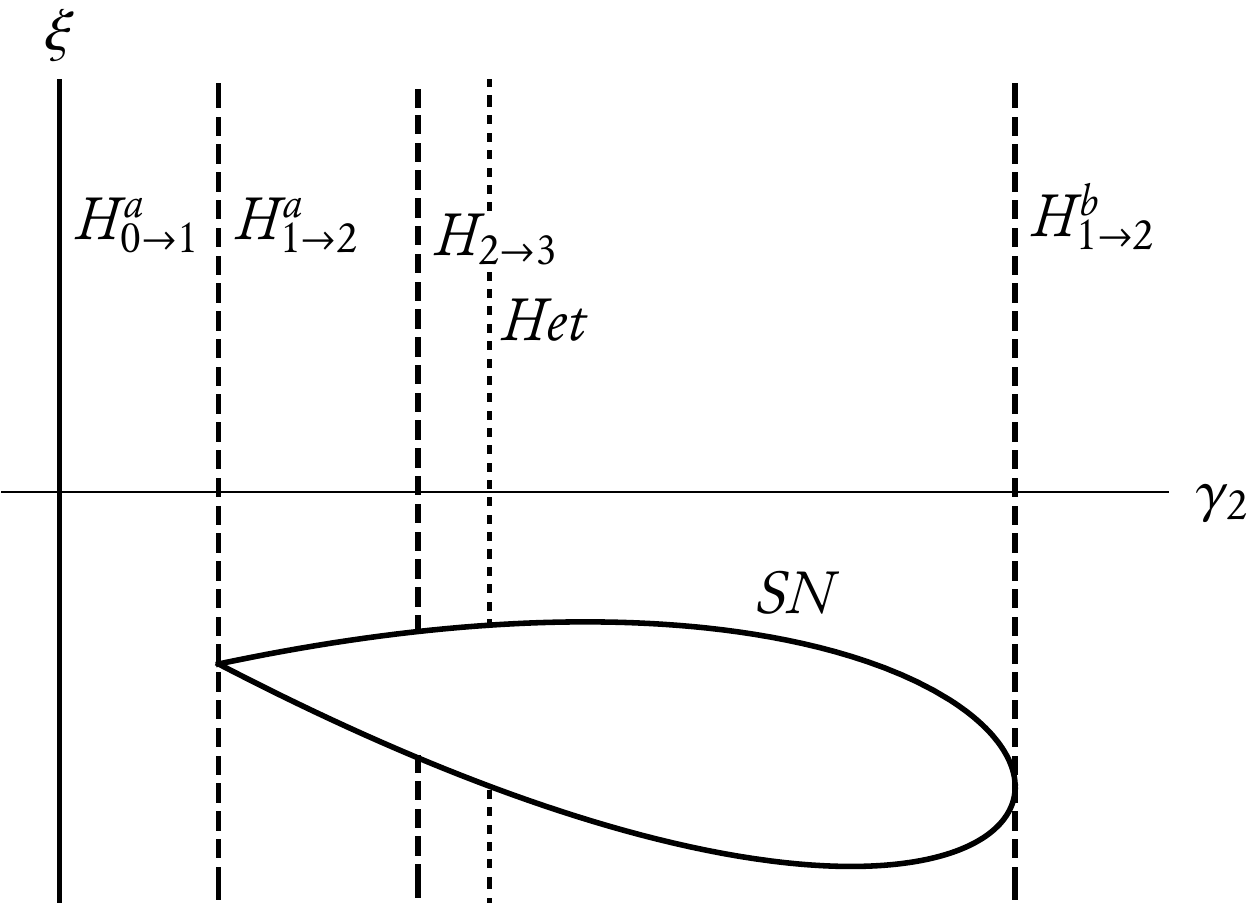}
    \subcaption{$1{:}2$ resonance}
  \end{minipage}
  \hfil
  \begin{minipage}[b]{0.425\linewidth}
    \centering
    \includegraphics[width=\textwidth]{%
      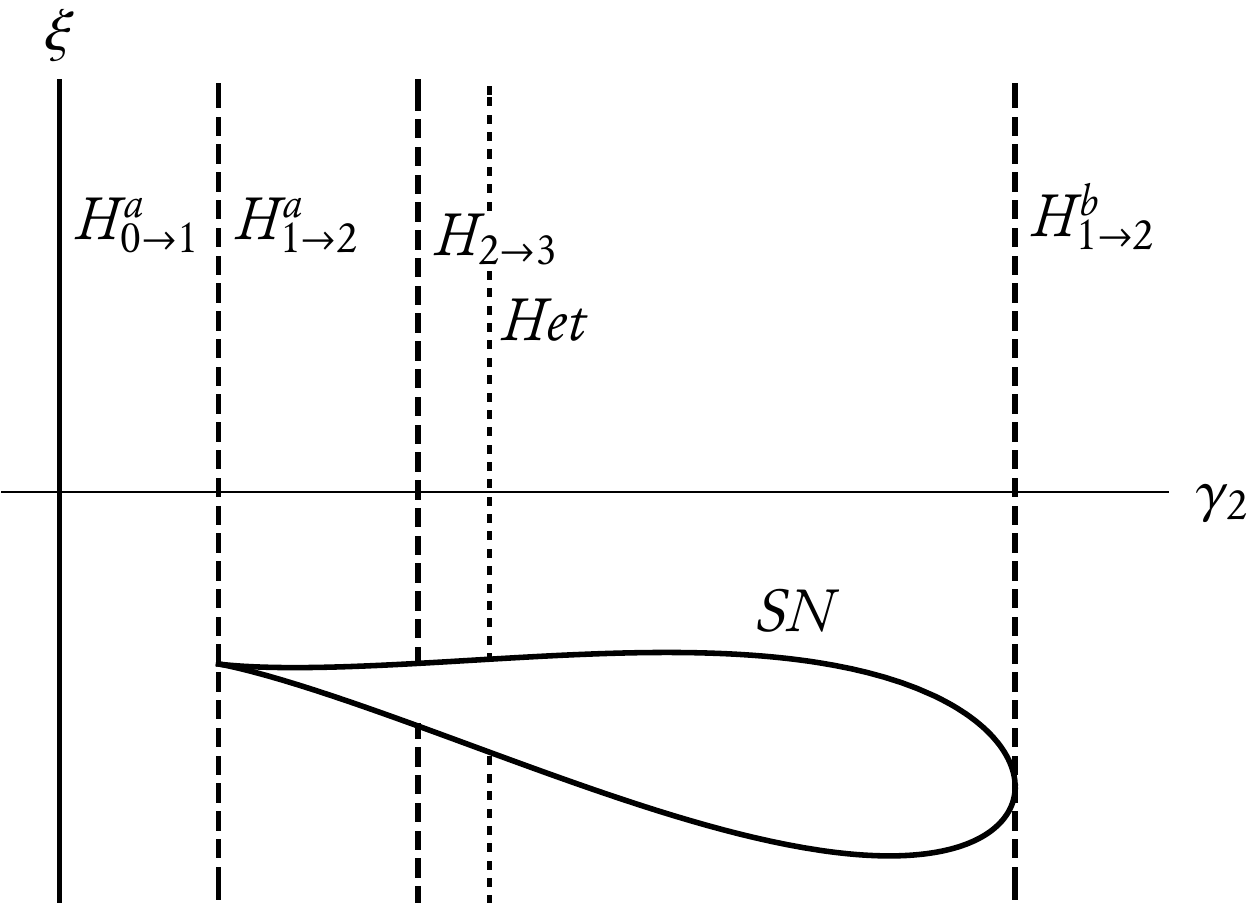}
    \subcaption{$1{:}3$ resonance} 
  \end{minipage}

  \begin{minipage}[b]{0.425\linewidth}
    \centering
    \includegraphics[width=\textwidth]{%
      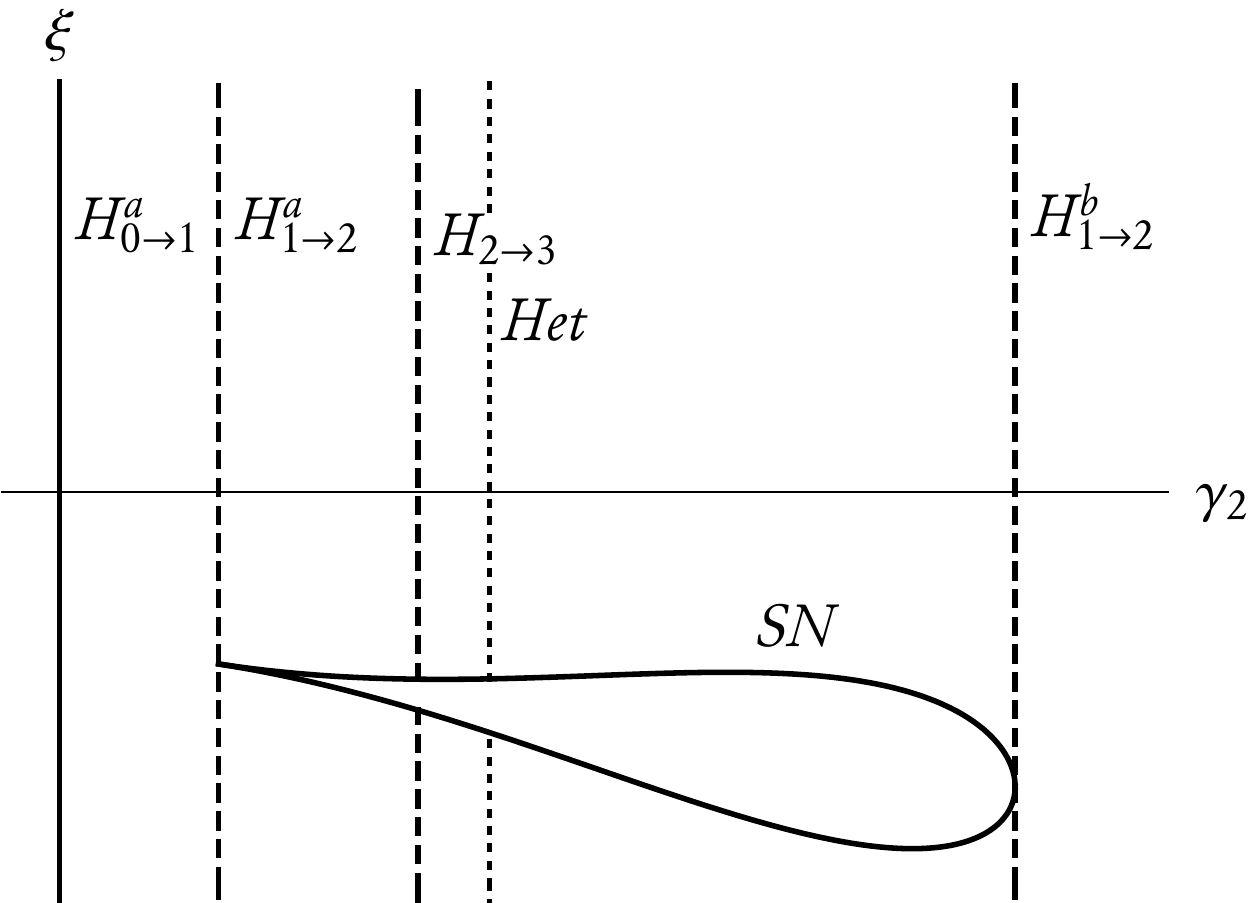}
    \subcaption{$1{:}4$ resonance}
  \end{minipage}
  \hfil
  \begin{minipage}[b]{0.425\linewidth}
    \centering
    \includegraphics[width=\textwidth]{%
      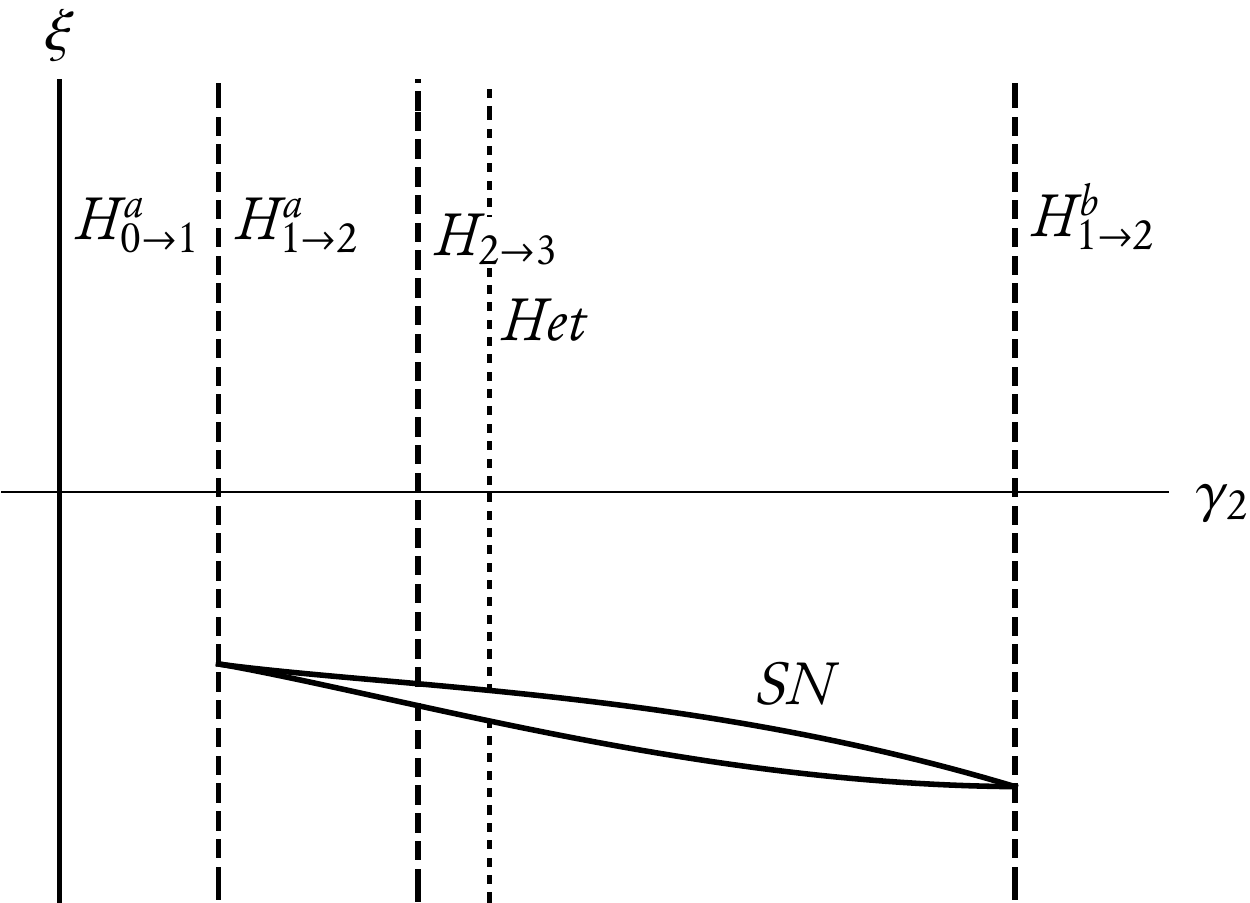}
    \subcaption{$2{:}3$ resonance} 
  \end{minipage}
  \caption{Vertical sections of the bifurcation diagram given in
    figure~\ref{fig:bd}.
    So we fix $\gamma_1$ and only consider the parameters $\gamma_2$
    and~$\xi$.
    Shape of the resonance droplets for the $1{:}2$, the $1{:}3$, the
    $1{:}4$ and the $2{:}3$~resonances.} 
  \label{fig:SN1}
\end{figure}

\samethought
Writing
$R^2 = \bar{\sigma}_1^{\ell_2} \bar{\sigma}_2^{\ell_1}/G_\ell$
and $m = (B - P A)^T \ell^{\perp}$, the condition that the line
$\inprod{m}{\psi} + \zeta = 0$ touches the circle $|\psi| = R$
implies that $\psi = \alpha m$ for a real scalar~$\alpha$ that
satisfies $|\alpha| = R/|m|$. 
For $\eps = 0$ the saddle-node bifurcation condition can then
be written as
\begin{equation}
  \label{eq:tangency}
  \zeta \;\; = \;\; \pm |m|R \;\; = \;\;
  \pm C \bar{\sigma}_1^{\ell_2/2} \bar{\sigma}_2^{\ell_1/2}.
\end{equation}
Since a line is flat and a circle has constant curvature, it is
clear that --- for sufficiently small~$\varepsilon > 0$ --- the
saddle-node bifurcation is nondegenerate. 

\begin{figure}[htp]
  \centering
  \begin{minipage}[b]{0.425\linewidth}
    \centering
    \includegraphics[width=\textwidth]{%
      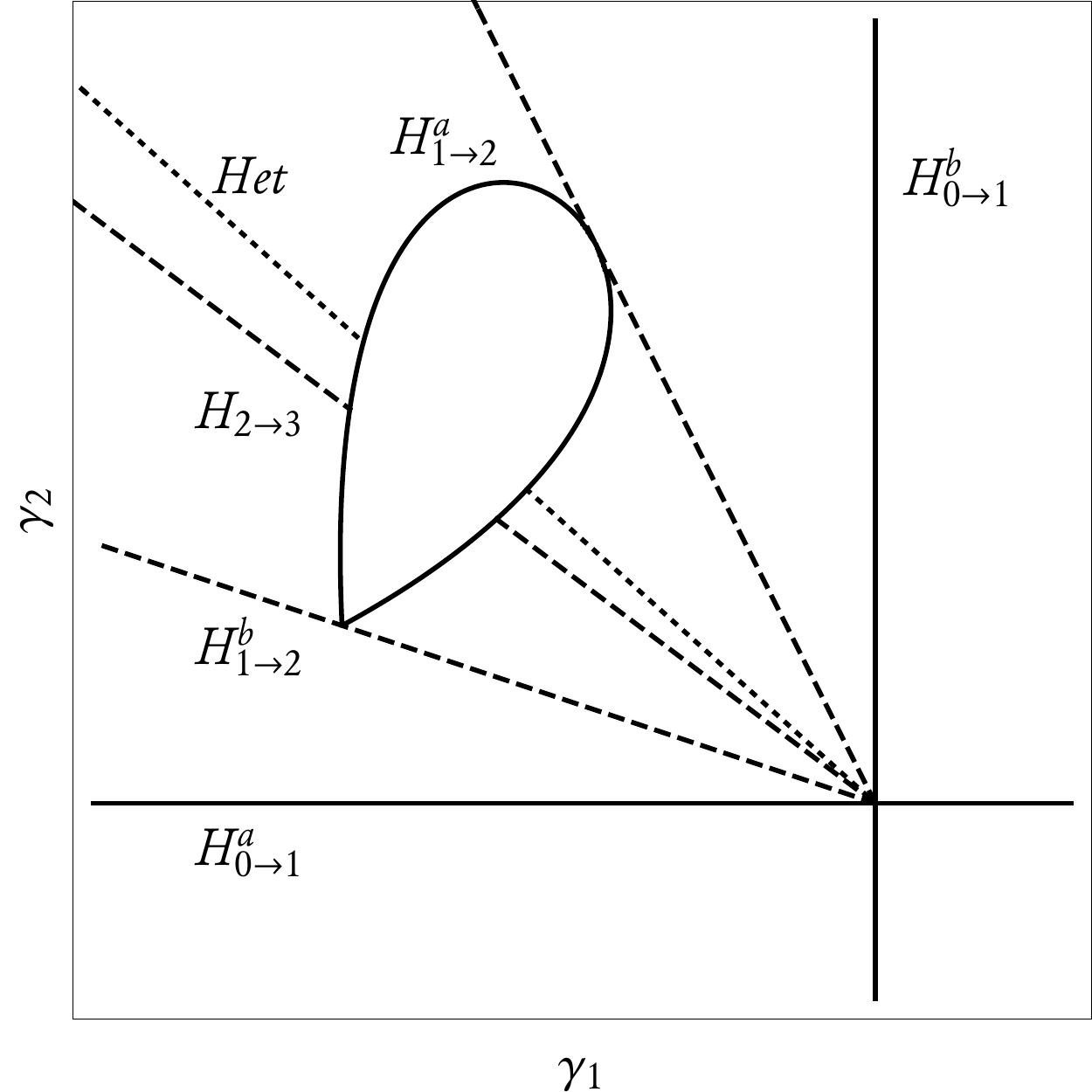}
    \subcaption{$1{:}2$ resonance}
  \end{minipage}
  \hfil
  \begin{minipage}[b]{0.425\linewidth}
    \centering
    \includegraphics[width=\textwidth]{%
      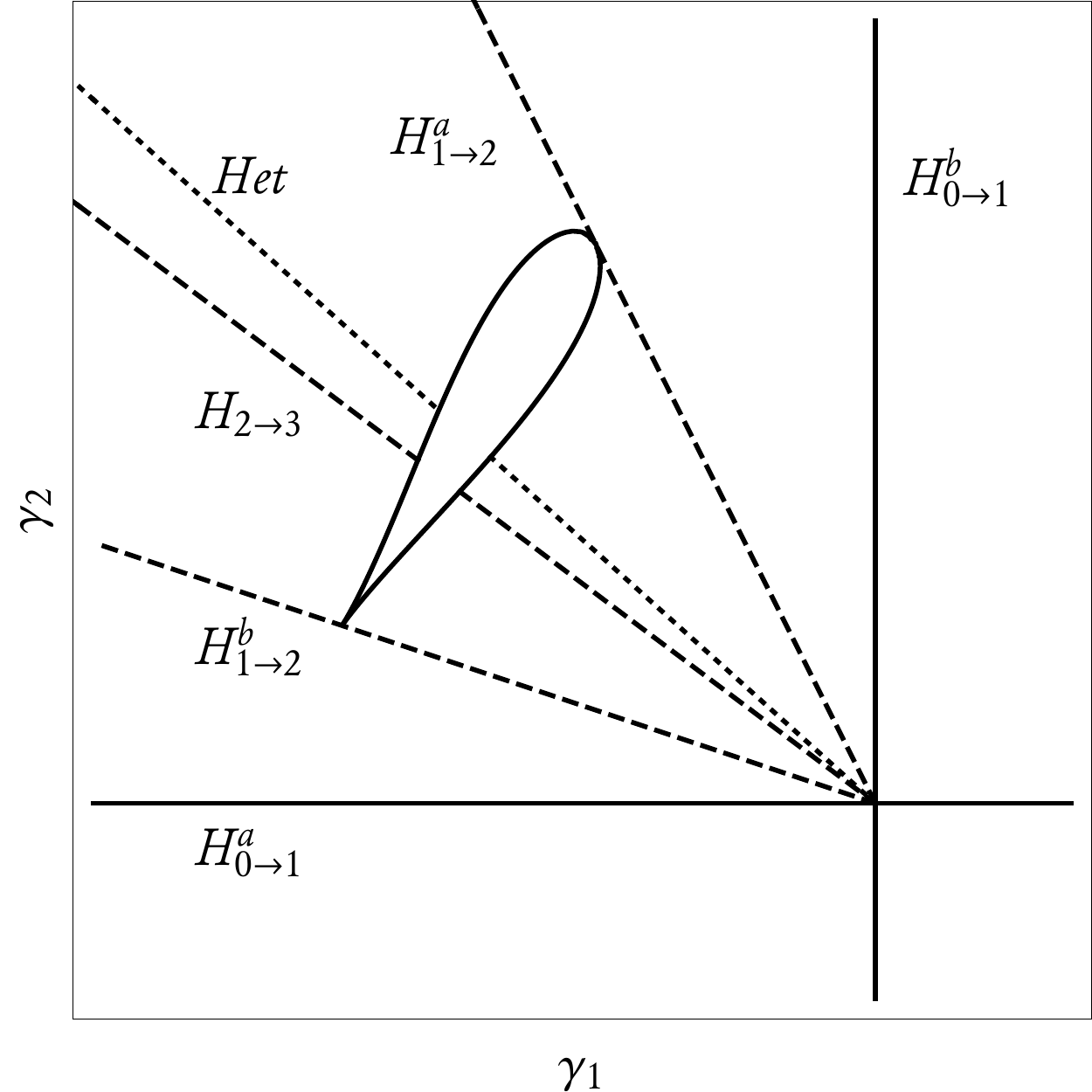}
    \subcaption{$1{:}3$ resonance} 
  \end{minipage}
  
  \vspace{2em}

  \begin{minipage}[b]{0.425\linewidth}
    \centering
    \includegraphics[width=\textwidth]{%
      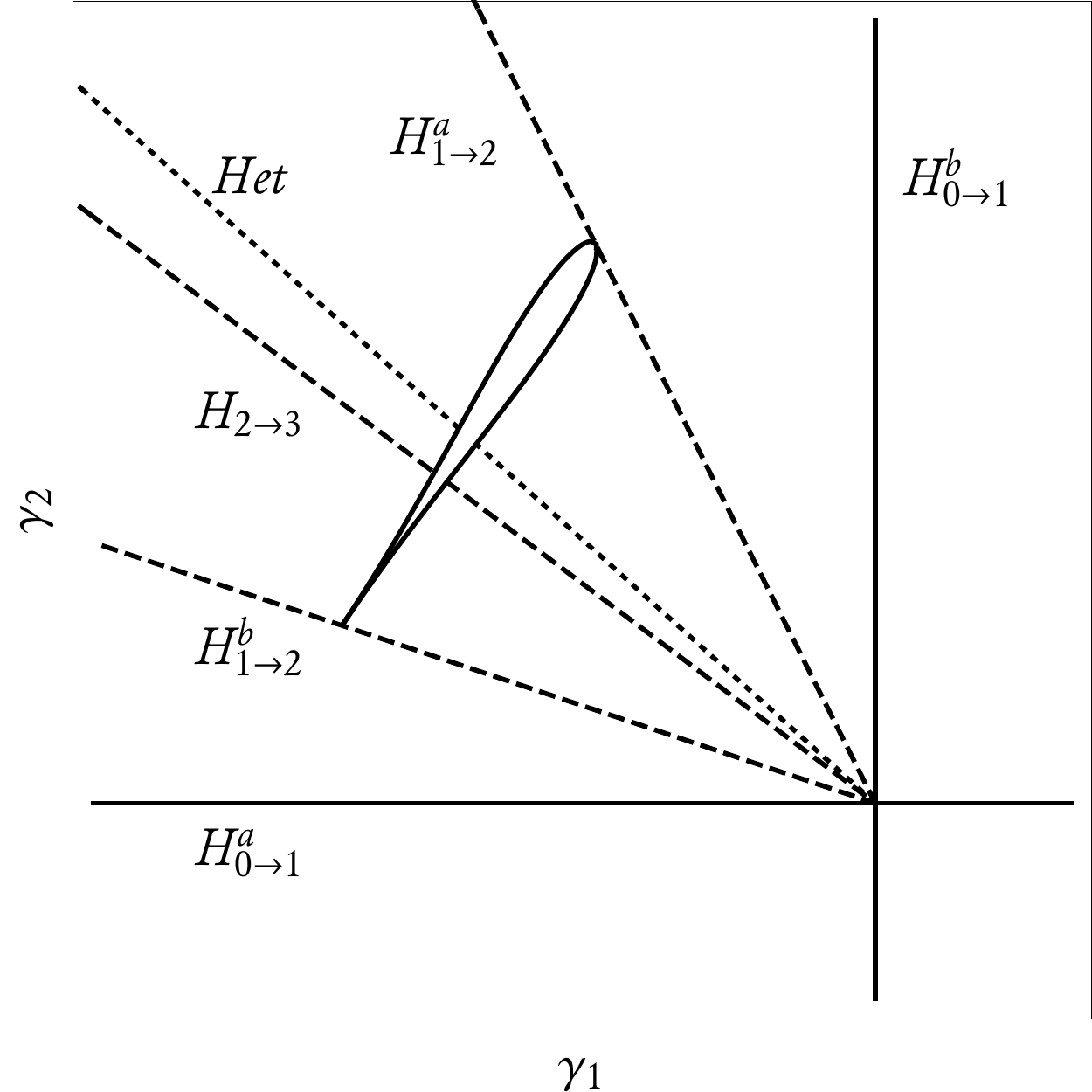}
    \subcaption{$1{:}4$ resonance}
  \end{minipage}
  \hfil
  \begin{minipage}[b]{0.425\linewidth}
    \centering
    \includegraphics[width=\textwidth]{%
      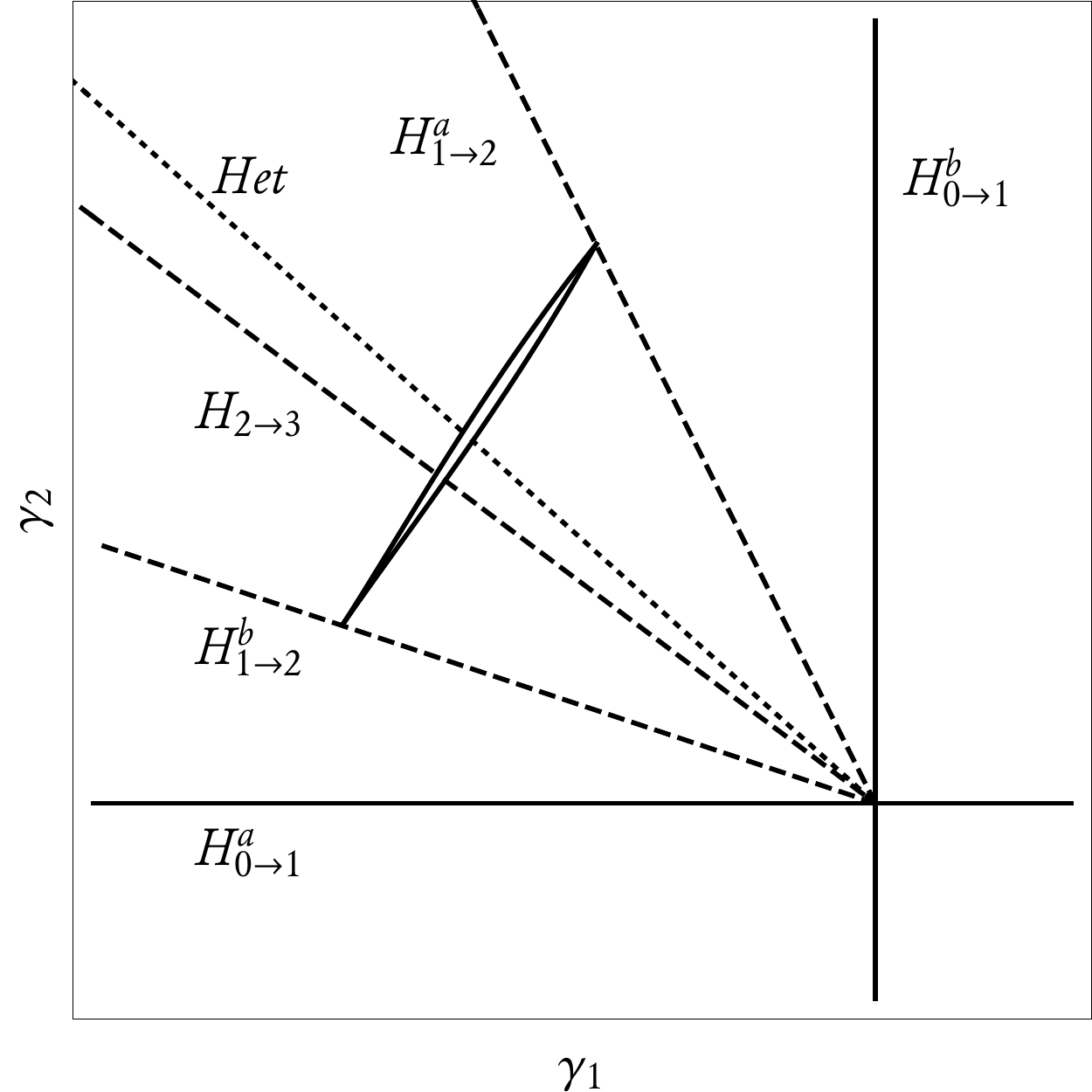}
    \subcaption{$2{:}3$ resonance} 
  \end{minipage}
  \caption{Horizontal sections of the bifurcation diagram given in
    figure~\ref{fig:bd}.
    So we fix $\xi$ and only consider the parameters $\gamma_1$
    and~$\gamma_2$.
    Shape of the resonance droplets for the $1{:}2$, the $1{:}3$, the
    $1{:}4$ and the $2{:}3$~resonances.} 
  \label{fig:SN2}
\end{figure}

\samethought
Equation~\eqref{eq:tangency} shows that the saddle-node
curves are the boundaries of a resonance droplet with width
proportional to $\bar{\sigma}_1^{\ell_2/2}$ close to the Hopf
curve $\bar{\sigma}_1 = 0$ and width proportional to
$\bar{\sigma}_2^{\ell_1/2}$ close to the Hopf curve
$\bar{\sigma}_2 = 0$. 
This is illustrated in figures \ref{fig:SN1} and~\ref{fig:SN2}.

\subsection{Putting it all together}
When the saddle-node bifurcation on the invariant circle
coincides with the occurrence of two unit Floquet multipliers of
the limit circle, we expect a fold-Hopf bifurcation to occur.
In figures \ref{fig:SN1} and~\ref{fig:SN2} these are the two
points where the line~$H_{2 \to 3}$ --- which ceases to exist
inside the resonance droplet --- meets the resonance droplet at the
boundary.

\section{Dynamics of the full system}
\label{sec:dynamicsofthefullsystem}

In the previous section we analysed the dynamics of the truncated
normal form on~$\R^4$, i.e.\ after reduction of the $\T^n$~symmetry.
Reconstructing the dynamics of the truncated normal form on
$\T^n \times \R^4$ amounts to restoring the driving
\begin{displaymath}
  \dot{x} \;\; = \;\; \omega(\mu)
            \; + \;
            \hat f(\tau_1, \tau_2,\mu)
            \; + \;
            a_0(\mu) \tau_3
            \; + \;
            c_0(\mu) \tau_4,
\end{displaymath}
which for $\tau = 0$ is conditionally periodic with frequency
vector~$\omega(\mu)$.
On~$\R^4$ we still have a $\T$~symmetry in the resonant case and
a $\T^2$~symmetry in the non-resonant case (or if we choose an
order of truncation lower than the order~$|\ell|$ of $\tau_3$
and~$\tau_4$).
Reducing the former takes us to the reduced phase space~$\sP$ and
reducing the latter takes us directly to the basis~$\sB$ of the
circle bundle~$\sP$.

\newthought
For general survey in the relevant \textsc{kam} theory, see
\cite{BHS96,BT09,BS10,CLB05}. 
The general philosophy is that we first are given a generic integrable
(i.e.~torus symmetric) system, where the dimension of the torus equals
$n$, $n+1$ or $n+2$.
  The Diophantine conditions \eqref{eq:diophantine} determine a
  `Cantorised' sub-bundle which is nowhere dense, but of large
  Lebesgue measure for small values of the gap parameter~$\Gamma$.
  Restricted to this set a Whitney smooth conjugation exists with a
  subset of the perturbed system.
  This result is called \textit{quasi-periodic stability}: structural
  stability restricted to a bundle of Diophantine quasi-periodic tori.
  The main question here is what happens in the gaps of this
  Cantorised bundle.
  In the present dissipative context hyperbolicity plays
  an important role in the perturbation analysis.

\subsection{No resonance}

The truncated normal form is not only independent of~$x$ but also
independent of $\tau_3$ and~$\tau_4$ if there is no
resonance~\eqref{eq:normal:resonance} with $\ell \neq 0$ between
the normal frequencies $\alpha_1(0)$ and~$\alpha_2(0)$, or if we
chose an order of truncation that is lower than the order~$|\ell|$
of $\tau_3$ and~$\tau_4$ in case there is a resonance $\ell \neq 0$
of the form~\eqref{eq:normal:resonance}.
Reconstructing the dynamics on~$\R^4$ from the reduced dynamics
on~$\sB$ is already implicit in figure~\ref{fig:hopf}.
Indeed, it is the reconstructed toral dynamics that turns the
pitchfork bifurcations that are actually visible in that figure
into the Hopf bifurcations $H^{a, b}_{0 \to 1}$
and~$H^{a, b}_{1 \to 2}$.
On $\T^n \times \R^4$ these then reconstruct to Hopf bifurcations
of $n$--tori and $(n+1)$--tori, respectively, resulting in tori of
dimensions $n+1$ and $n+2$.

\newthought
For $n = 0$ persistence of the latter --- when perturbing from the
normal form back to the original system --- has been proved
in~\cite{li16}, as has been persistence of the $3$--tori resulting
from the quasi-periodic Hopf bifurcation~$H_{2 \to 3}$.
Persistence of the three types of Hopf bifurcations of $n$--tori,
$(n+1)$--tori and $(n+2)$--tori also follows
from~\cite{BBH90, CLB05}, see furthermore~\cite{BHW20}.
We remark that the specialised approach in~\cite{li16} yields the
same order $\varepsilon^{1/18}$ as the general treatment
in~\cite{BBH90}. 
 Indeed, in~\cite{BBH90} we have a general normal form as in
 Theorem~\ref{thm:normalform},
 where we have to go up to $8$th order.
 Accounting each normalisation step by $\sqrt{\varepsilon}$,
 the remainder terms turn out to be of order~$\varepsilon^{1/18}$. 

\samethought
The gaps left open when Cantorising the attracting tori to prove
persistence using \KAM~theory can be closed using normal
hyperbolicity to obtain invariant tori on which the flow is no
longer quasi-periodic, compare with~\cite{BHW18} and results cited
therein.
There seem to be no claims concerning the heteroclinic bifurcation
in the literature.

\subsection{Higher order resonances}

Here we normalise op to order $|\ell| \geq 5$.
This makes the $1{:}4$ and $2{:}3$~resonances still special among
the higher order resonances as normalising up to $5$th order to
resolve the Hopf bifurcation~$H_{2 \to 3}$ automatically includes
$\tau_3$ and~$\tau_4$ into the normal form.
Note that the opening of the resonance bubbles is dictated by
$\ell_1$ and~$\ell_2$ separately, see e.g.\ figures \ref{fig:SN1}c
and~\ref{fig:SN2}c.

\newthought
Where the line~$H_{2 \to 3}$ meets the resonance bubble we have a
fold-Hopf bifurcation; in fact the line ceases to exist inside the
resonance bubble.
Next to the Hopf bifurcations~$H^{a, b}_{0 \to 1}$ of $n$--tori,
for the persistence of which we again refer to
\cite{BBH90, CLB05, BHW20}, this leaves us with six special points:
the two Hopf bifurcations~$H^{a, b}_{1 \to 2}$ of $(n+1)$--tori with
a normal-internal resonance, the two fold-Hopf bifurcations and the
two heteroclinic bifurcations at the boundary of the resonance bubble
(which we do not further comment upon).

\subsection{The $1{:}3$~resonance}

In this case $|\ell| - 2 = 2$, which means that we have two competing
influences in the normal form.
Hence, the relative strengths of the coefficients $\tilde{p}$ and~$A\psi$
in \eqref{eq:sigma-psi}
decides whether we have dynamics similar to the
previous subsection or dynamics similar to the following subsection.

\subsection{The $1{:}2$~resonance}

Although the invariants $\tau_3$ and~$\tau_4$ have order $|\ell| = 3$
in this case, we still work with a $5$th order normal form to resolve
the Hopf bifurcation~$H_{2 \to 3}$.
In this way not only the linear terms $\tau_3$ and~$\tau_4$ enter the
normal form, but also the four $5$th order terms $\tau_i \tau_j$ with
$i = 1, 2$ and $j = 3, 4$.
We wonder whether the subordinate fold-Hopf bifurcation
requests even higher order terms of the normal form
or whether the $5$th order terms are already sufficient,
again see~\cite{BHW20}.

\newthought
At the end points of the resonance bubble we have periodic Hopf
bifurcations~$H^{a, b}_{1 \to 2}$ with normal-internal resonances
$1{:}1$ and~$1{:}2$, respectively.

\section{Final remarks}
\label{sec:conclusions}

The $\T^1$ or $\T^2$ symmetric normal form analysis
forms a kind of skeleton of the total dynamics, which
occurs when adding the non-symmetric, higher order terms that are
flat.

\samethought
First note that normal forms of any finite order keep the toroidal symmetry.
That means that flat terms are of infinite order and we expect these
terms to be exponentially small in the real analytic
context~\cite{BT89,BR01,BHW18}.
In the present dissipative context the higher order resonances
in the Diophantine conditions~\eqref{eq:diophantine}
give rise to infinite arrays of `small' resonance bubbles
\cite{bro03,BB87,BBH90,BHW18,BHW20,che85,che85a,che87}, 
where the array `runs along' the bifurcation manifolds as given by the
quasi-periodic sub-bundle, compare with the discussion at the
beginning of section~\ref{sec:dynamicsofthefullsystem}.
The paper~\cite{BHW18} exactly uses the quasi-periodic Hopf
bifurcation as a leading example.

\begin{figure}
  \centering
\includegraphics[height=0.3\textheight]{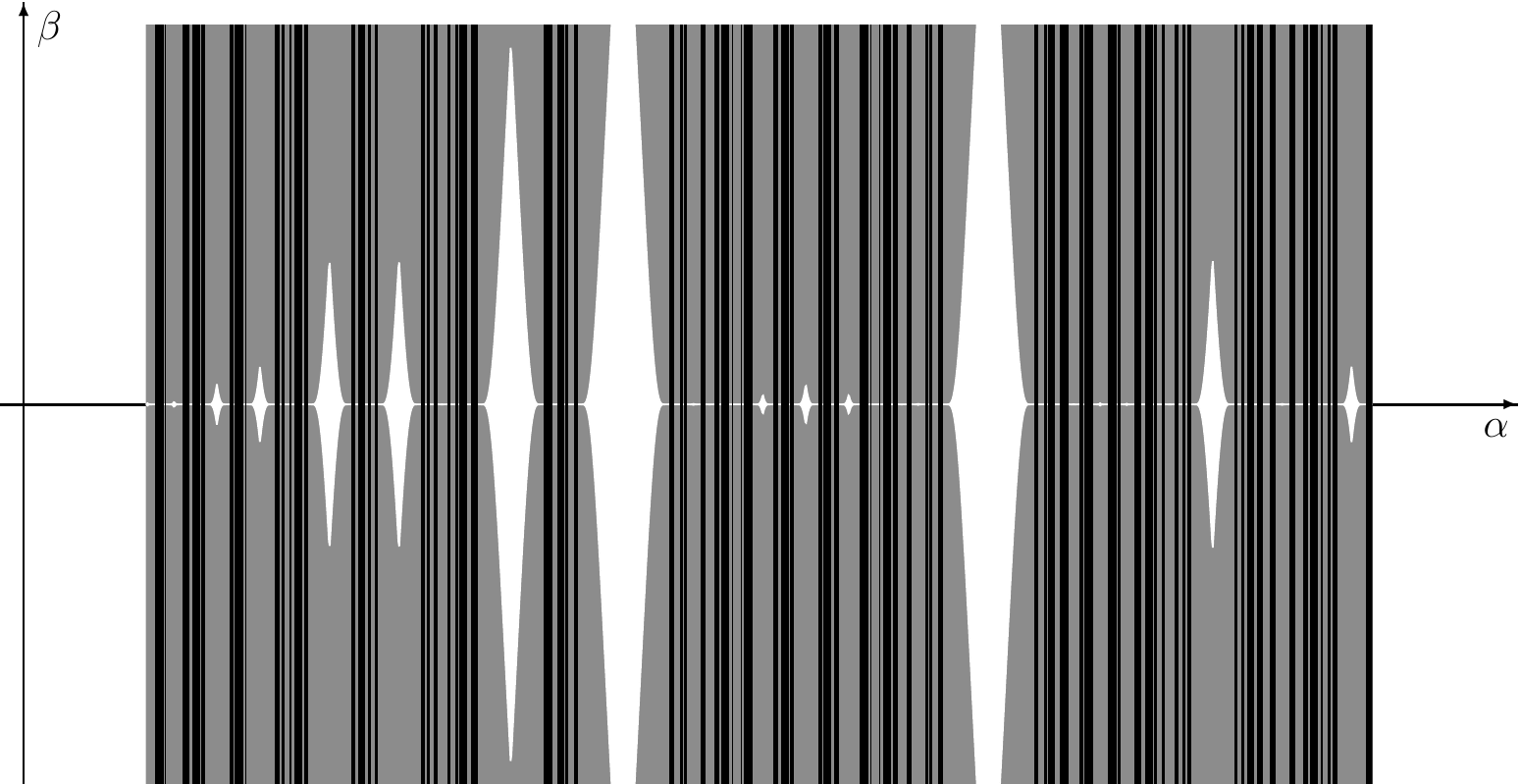}
\begin{small}
\caption{From~\cite{BHW18}. 
Parameter sets with persistent invariant tori and the subset of
persistent quasi-periodic invariant tori.
The latter (black) set is the product of a Cantor set in the $\alpha$--direction
and the real line in the $\beta$--direction: it is nowhere dense.
The (white) complement of the former set (which includes the resonance set
on the $\beta$--axis) is a ‘set of ignorance’.
\label{fig:qph}}
\end{small}
\end{figure}

\newthought
Inside the bubbles the complexity is large both in the
  dissipative and conservative context, which is illustrated
  by the papers~\cite{BSV08a,BSV08b,BSV10}. 
  






\end{document}